# Testing Endogeneity with High Dimensional Covariates[*]

Zijian Guo[1], Hyunseung Kang[2], T. Tony Cai[3], and Dylan S. Small[3]

[1]Department of Statistics and Biostatistics, Rutgers University
[2]Department of Statistics, University of Wisconsin-Madison
[3]Department of Statistics, The Wharton School, University of Pennsylvania

## Abstract

Modern, high dimensional data has renewed investigation on instrumental variables (IV) analysis, primarily focusing on estimation of effects of endogenous variables and putting little attention towards specification tests. This paper studies in high dimensions the Durbin-Wu-Hausman (DWH) test, a popular specification test for endogeneity in IV regression. We show, surprisingly, that the DWH test maintains its size in high dimensions, but at an expense of power. We propose a new test that remedies this issue and has better power than the DWH test. Simulation studies reveal that our test achieves near-oracle performance to detect endogeneity.

---

[*]Address for correspondence: Zijian Guo, Department of Statistics and Biostatistics, Rutgers University, USA. Phone: (848)445-2690. Fax: (732)445-3428. Email: zijguo@stat.rutgers.edu.

# 1 Introduction

## 1.1 Endogeneity Testing with High Dimensional Data

Recent growth in both the size and dimension of the data has led to a resurgence in analyzing instrumental variables (IV) regression in high dimensional settings (Belloni et al., 2012, 2013, 2011a; Chernozhukov et al., 2014, 2015; Fan and Liao, 2014; Gautier and Tsybakov, 2011) where the number of regression parameters, especially those associated with exogenous covariates, is growing with, and may exceed, the sample size.[1] The primary focus in these works has been providing tools for estimation and inference of a single endogenous variable's effect on the outcome under some low-dimensional structural assumptions on the structural parameters associated with the instruments and the covariates, such as sparsity. (Belloni et al., 2012, 2013, 2011a; Chernozhukov et al., 2014, 2015; Gautier and Tsybakov, 2011). This line of work has generally not focused on specification tests in the high dimensional IV setting.

The main goal of this paper is to study the high dimensional behavior of one of the most common specification tests in IV regression, the test for endogeneity, which assumes the validity of the IV and tests whether the included endogenous variable (e.g., a treatment variable) is actually exogenous. Historically, the most widely used test for endogeneity is the Durbin-Wu-Hausman test (Durbin, 1954; Hausman, 1978; Wu, 1973), hereafter called the DWH test, and is widely implemented in software, such as ivreg2 in Stata (Baum et al., 2007). The DWH test detects the presence of endogeneity in the structural model by studying the difference between the ordinary least squares (OLS) estimate of the structural parameters in the IV regression to that of the two-stage least squares (TSLS) under the null hypothesis of no endogeneity; see Section 2.3 for the exact characterization of the DWH test. In low dimensional settings, the primary requirements for the DWH test to correctly control Type I error are having instruments that are (i) strongly associated with the included endogenous variable, often called strong instruments, and (ii) exogenous to the structural errors[2], often referred to as valid instruments (Murray, 2006). When instruments are not strong, Staiger and Stock (1997) showed that the DWH test that used the TSLS estimator for variance, developed by Durbin (1954) and Wu (1973), had distorted size under the null hypothesis while the DWH test that used the OLS estimator for variance, developed by Hausman (1978), had proper size. When instruments are invalid, which is perhaps a bigger concern in practice (Conley et al., 2012; Murray, 2006), the DWH test will usually fail because the TSLS estimator is inconsistent under the null hypothesis; see the Supplementary materials for a simple theoretical justification of this phenomenon. Indeed, some recent work with high dimensional data (Belloni et al., 2012; Chernozhukov et al., 2015) advocated conditioning on many, possibly high dimensional, exogenous covariates to make instruments more plausibly valid.[3] However, while adding additional covariates can potentially make instruments more plausibly valid, it is unclear what price one has to pay with respect to the performance of specification tests like the DWH test.

---

[1] In the paper, we use the term "high dimensional setting" more broadly where the number of parameters is growing with the sample size; see Sections 3 and 4.3 for details and examples. Note that the modern usage of the term "high dimensional setting" where the sample size exceeds the parameter is one case of this broader setting.

## 1.2 Prior Work and Contribution

Prior work in analyzing the DWH test in instrumental variables is diverse. Estimation and inference under weak and/or many instruments are well documented (Andrews et al., 2007; Bekker, 1994; Bound et al., 1995; Chao and Swanson, 2005; Dufour, 1997; Han and Phillips, 2006; Hansen et al., 2008; Kleibergen, 2002; Moreira, 2003; Morimune, 1983; Nelson and Startz, 1990; Newey and Windmeijer, 2005; Staiger and Stock, 1997; Stock and Yogo, 2005; Wang and Zivot, 1998; Zivot et al., 1998). In particular, when the instruments are weak, the behavior of the DWH test under the null depends on the variance estimate (Doko Tchatoka, 2015; Nakamura and Nakamura, 1981; Staiger and Stock, 1997). Other works study the behavior of the DWH test under different strengths of instruments and/or weak instrument asymptotics (Hahn et al., 2011; Staiger and Stock, 1997) and under a two-stage testing scheme (Guggenberger, 2010). Some recent work extended the specification test to handle growing number of instruments (Chao et al., 2014; Hahn and Hausman, 2002; Lee and Okui, 2012). Other recent works extended specification tests based on overidentification (Hahn and Hausman, 2005; Hausman et al., 2005) and to heteroskedastic data (Chmelarova et al., 2007). Fan et al. (2015) considered testing endogeneity in the high dimensional non-IV setting and approximated the null distribution of their test statistic by the bootstrap; the distribution under the alternative was not identifiable. None of these works have characterized the properties of the DWH test used in IV regression under the high dimensional setting.

---

[2]The term exogeneity is sometimes used in the IV literature to encompass two assumptions, (a) independence of the IVs to the disturbances in the structural model and (b) IVs having no direct effect on the outcome, sometimes referred to as the exclusion restriction (Angrist et al., 1996; Holland, 1988; Imbens and Angrist, 1994). As such, an instrument that is perfectly randomized from a randomized experiment may not be exogenous in the sense that while the instrument is independent to any structural error terms, the instrument may still have a direct effect on the outcome.

[3]For example, in Section 7 of the empirical example of Belloni et al. (2012), the authors studied the effect of federal appellate court decisions on economic outcomes by using the random assignment of judges to decide appellate cases. They state that once the distribution of characteristics of federal circuit court judges in a given circuit-year is controlled for, "the realized characteristics of the randomly assigned three-judge panel should be unrelated to other factors besides judicial decisions that may be related to economic outcomes" (page 2405). More broadly, in empirical practice, adding covariates to make IVs more plausibly valid is commonplace; see Card (1999), Cawley et al. (2013), and Kosec (2014) for examples as well as review papers in epidemiology and causal inference by Hernán and Robins (2006) and Baiocchi et al. (2014).

Our main contributions are two-fold. First, we characterize the behavior of the popular DWH test in high dimensions. The theoretical analysis reveals that the DWH test actually controls Type I error at the correct level in high dimensions, but pays a significant price with respect to power, especially for small to moderate degrees of endogeneity; we also confirm our finding numerically with a simulation study of the empirical power of the DWH test. Our finding also suggests that, although conditioning on a large number of covariates makes instruments more plausibly valid, the power of the DWH test is reduced because of the large number of covariates. Second, we remedy the low power of the DWH test by presenting a simple and improved endogeneity test that is robust to high dimensional covariates and/or instruments and that works in settings where the number of structural parameters is allowed to exceed the sample size. In particular, our new endogeneity test applies a hard thresholding procedure to popular estimators for reduced-form models, such as OLS in low dimensions or bias-corrected Lasso estimators in high dimensions (see Section 4.1 for details). This hard thresholding procedure is an essential step of the new endogeneity test, where relevant instruments are selected for testing endogeneity. We also highlight that the success of the proposed endogeneity test does not require the correct selection of all relevant instruments. That is, even if the relevant instruments are not correctly selected, the proposed testing procedure still controls Type I error and achieves non-trivial power under regularity conditions. Additionally, we briefly discuss an extension of our endogeneity test to incorporate invalid instruments, especially when many covariates are conditioned upon to avoid invalid IVs.

This paper is closely connected to the paper Guo et al. (2016) by the same group of authors, where Guo et al. (2016) proposed confidence intervals for the treatment effect in the presence of both high-dimensional instruments, covariates and invalid instruments. The current paper considers a related but different problem about endogeneity testing and extends the idea proposed in Guo et al. (2016) to testing endogeneity in high dimensional settings. In particular, we are the first to provide a test for endogeneity when $n < p$. In addition, the characterization of the power of DWH test in high dimensions and the technical tools used in the paper are new. In particular, the technical tools can be used to study other specification tests like the Sargan or the J test (e.g. Sargan test (Hansen, 1982; Sargan, 1958)) in high dimensions.

We conduct simulation studies comparing the performance of our new test with the usual DWH test and apply the proposed endogeneity test to an empirical data analysis following Belloni et al. (2012, 2014). We find that our test has the desired size and has better power than the DWH test for all degrees of endogeneity and performs similarly to the oracle DWH test which knows the support of relevant instruments and covariates a priori. In the supplementary materials, we also present technical proofs and extended simulation studies that further examine the power of our test.

## 2 Instrumental Variables Regression and the DWH Test

### 2.1 Notation

For any vector $v \in \mathbb{R}^p$, $v_j$ denotes its $j$th element, and $\|v\|_1$, $\|v\|_2$, and $\|v\|_\infty$ denote the $1, 2$ and $\infty$-norms, respectively. Let $\|v\|_0$ denote the number of non-zero elements in $v$ and define $\mathrm{supp}(v) = \{j : v_j \neq 0\} \subseteq \{1, \ldots, p\}$. For any $n \times p$ matrix $M$, denote the $(i, j)$ entry by $M_{ij}$, the $i$th row by $M_{i.}$, the $j$th column by $M_{.j}$, and the transpose of $M$ by $M'$. Also, given any $n \times p$ matrix $M$ with sets $I \subseteq \{1,\ldots,n\}$ and $J \subseteq \{1,\ldots,p\}$ denote $M_{IJ}$ as the submatrix of $M$ consisting of rows specified by the set $I$ and columns specified by the set $J$, $M_{I.}$ as the submatrix of $M$ consisting of rows indexed by the set $I$ and all columns, and $M_{.J}$ as the submatrix of $M$ consisting of columns specified by the set $J$ and all rows. Also, for any $n \times p$ full-rank matrix $M$, define the orthogonal projection matrices $P_M = M(M'M)^{-1}M'$ and $P_{M^\perp} = I - M(M'M)^{-1}M'$ where $P_M + P_{M^\perp} = I$ and I is an identity matrix. For a $p \times p$ matrix $\Lambda$, $\Lambda \succ 0$ denotes that $\Lambda$ is a positive definite matrix. For any $p \times p$ positive definite $\Lambda$ and set $J \subseteq \{1,\ldots,p\}$, let $\Lambda_{J|J^C} = \Lambda_{JJ} - \Lambda_{JJ^C}\Lambda_{J^C J^C}^{-1}\Lambda_{J^C J}$ denote the submatrix $\Lambda_{JJ}$ adjusted for the columns in the complement of the set $J$, $J^C$.

For a sequence of random variables $X_n$ indexed by $n$, we use $X_n \xrightarrow{p} X$ to represent that $X_n$ converges to $X$ in probability. For a sequence of random variables $X_n$ and numbers $a_n$, we define $X_n = o_p(a_n)$ if $X_n/a_n$ converges to zero in probability and $X_n = O_p(a_n)$ if for every $c_0 > 0$, there exists a finite constant $C_0$ such that $\mathbf{P}(|X_n/a_n| \geq C_0) \leq c_0$. For any two sequences of numbers $a_n$ and $b_n$, we will write $b_n \ll a_n$ if $\limsup b_n/a_n = 0$.

For notational convenience, for any $\alpha$, $0 < \alpha < 1$, let $\Phi$ and $z_{\alpha/2}$ denote, respectively, the cumulative distribution function and $\alpha/2$ quantile of a standard normal distribution. Also, for any $B \in \mathbb{R}$, we define the function $G(\alpha, B)$ to be the tail probabilities of a normal distribution shifted by $B$, i.e.

$$G(\alpha, B) = 1 - \Phi(z_{\alpha/2} - B) + \Phi(-z_{\alpha/2} - B). \tag{1}$$

We use $\chi^2_\alpha(d)$ to denote the $1 - \alpha$ quantile of the Chi-squared distribution with $d$ degrees of freedom.

## 2.2 Model and Definitions

Suppose we have $n$ individuals where for each individual $i = 1,\ldots,n$, we measure the outcome $Y_i$, the included endogenous variable $D_i$, $p_z$ candidate instruments $Z'_{i.}$, and $p_x$ exogenous covariates $X'_{i.}$ in an i.i.d. fashion. We denote $W'_{i.}$ to be concatenated vector of $Z'_{i\cdot}$ and $X'_{i\cdot}$ with dimension $p = p_z + p_x$. The columns of the matrix $W$ are indexed by two sets, the set $\mathcal{I} = \{1, \ldots, p_z\}$, which consists of all the $p_z$ candidate instruments, and the set $\mathcal{I}^C = \{p_z + 1, \ldots, p\}$, which consists of the $p_x$ covariates. The variables $(Y_i, D_i, Z_i, X_i)$ are governed by the following structural model.

$$Y_i = D_i\beta + X'_{i\cdot}\phi + \delta_i, \qquad E(\delta_i \mid Z_{i\cdot}, X_{i\cdot}) = 0 \qquad (2)$$

$$D_i = Z'_{i\cdot}\gamma + X'_{i\cdot}\psi + \epsilon_i, \qquad E(\epsilon_i \mid Z_{i\cdot}, X_{i\cdot}) = 0 \qquad (3)$$

where $\beta, \phi, \gamma$, and $\psi$ are unknown parameters in the model and without loss of generality, we assume the variables are centered to mean zero.[4] Let the population covariance matrix of $(\delta_i, \epsilon_i)$ be $\Sigma$, with $\Sigma_{11} = \text{Var}(\delta_i \mid Z_{i\cdot}, X_{i\cdot})$, $\Sigma_{22} = \text{Var}(\epsilon_i \mid Z_{i\cdot}, X_{i\cdot})$, and $\Sigma_{12} = \Sigma_{21} = \text{Cov}(\delta_i, \epsilon_i \mid Z_{i\cdot}, X_{i\cdot})$. Let the second order moments of $W_{i\cdot}$ be $\Lambda = E(W_{i\cdot}W'_{i\cdot})$ and let $\Lambda_{\mathcal{I}|\mathcal{I}^c}$ denote the adjusted covariance of variables belonging to the index set $\mathcal{I}$. Let $\omega$ represent all parameters $\omega = (\beta, \pi, \phi, \gamma, \psi, \Sigma)$ and define the parameter space

$$\Omega = \{\omega = (\beta, \pi, \phi, \gamma, \psi, \Sigma) : \beta \in \mathbb{R}, \pi, \gamma \in \mathbb{R}^{p_z}, \phi, \psi \in \mathbb{R}^{p_x}, \Sigma \in \mathbb{R}^{2\times 2}, \Sigma \succ 0\}. \qquad (4)$$

Finally, we denote $s_{x2} = \|\phi\|_0$, $s_{z1} = \|\gamma\|_0$, $s_{x1} = \|\psi\|_0$ and $s = \max\{s_{x2}, s_{z1}, s_{x1}\}$.

We also define relevant and irrelevant instruments. This is, in many ways, equivalent to the notion that the instruments $Z_{i\cdot}$ are associated with the endogenous variable $D_i$, except we use the support of a vector to define instruments' association to the endogenous variable; see Breusch et al. (1999); Hall and Peixe (2003), and Cheng and Liao (2015) for some examples in the literature of defining relevant and irrelevant instruments based on the support of a parameter.

**Definition 1.** *Suppose we have $p_z$ instruments along with model (3). We say that instrument $j = 1, \ldots, p_z$ is relevant if $\gamma_j \neq 0$ and irrelevant if $\gamma_j = 0$. Let $\mathcal{S} \subseteq \mathcal{I} = \{1, 2, \cdots, p_z\}$ denote the set of relevant IVs.*

Finally, for $\mathcal{S}$, the set of relevant IVs, we define the concentration parameter, a common

---

[4]Mean-centering is equivalent to adding a constant 1 term (i.e. intercept term) in $X'_{i\cdot}$; see Section 1.4 of Davidson and MacKinnon (1993) for details.

measure of instrument strength,

$$C(\mathcal{S}) = \frac{\gamma'_{\mathcal{S}} \Lambda_{\mathcal{S}|\mathcal{S}^C} \gamma_{\mathcal{S}}}{|\mathcal{S}|\Sigma_{22}}. \tag{5}$$

If all the instruments were relevant, then $\mathcal{S} = \mathcal{I}$ and equation (5) is the usual definition of concentration parameter in Bound et al. (1995); Mariano (1973); Staiger and Stock (1997) and Stock and Wright (2000) using population quantities, i.e. $\Lambda_{\mathcal{S}|\mathcal{S}^C}$. In particular, $C(\mathcal{S})$ corresponds exactly to the quantity $\lambda'\lambda/K_2$ on page 561 of Staiger and Stock (1997) when $n = 1$ and $K_1 = 0$. Without using population quantities, the function $nC(\mathcal{S})$ roughly corresponds to the usual concentration parameter using the sample estimator of $\Lambda_{\mathcal{S}|\mathcal{S}^C}$. However, if only a subset of all instruments are relevant so that $\mathcal{S} \subset \mathcal{I}$, then the concentration parameter in equation (5) represents the strength of instruments for that subset $\mathcal{S}$, adjusted for the exogenous variables in its complement $\mathcal{S}^C$. Regardless, like the usual concentration parameter, a high value of $C(\mathcal{S})$ represents strong instruments in the set $\mathcal{S}$ while a low value of $C(\mathcal{S})$ represents weak instruments.

## 2.3 The DWH Test

Consider the following hypotheses for detecting endogeneity in models (2) and (3),

$$H_0 : \Sigma_{12} = 0 \quad \text{versus} \quad H_1 : \Sigma_{12} \neq 0. \tag{6}$$

The DWH test tests the hypothesis of endogeneity in equation (6) by comparing two consistent estimators of $\beta$ under the null hypothesis $H_0$ (i.e. no endogeneity) with different efficiencies. Formally, the DWH test statistic, denoted as $Q_{\text{DWH}}$, is the quadratic difference between the OLS estimator of $\beta$, $\widehat{\beta}_{\text{OLS}} = (D'P_{X^\perp}D)^{-1}D'P_{X^\perp}Y$, and the TSLS estimator of $\beta$, $\widehat{\beta}_{\text{TSLS}} = (D'(P_W - P_X)D)^{-1}D'(P_W - P_X)Y$.

$$Q_{\text{DWH}} = \frac{(\widehat{\beta}_{\text{TSLS}} - \widehat{\beta}_{\text{OLS}})^2}{\widehat{\text{Var}}(\widehat{\beta}_{\text{TSLS}}) - \widehat{\text{Var}}(\widehat{\beta}_{\text{OLS}})}. \tag{7}$$

The terms $\widehat{\text{Var}}(\widehat{\beta}_{\text{OLS}})$ and $\widehat{\text{Var}}(\widehat{\beta}_{\text{TSLS}})$ are standard error estimates of the OLS and TSLS estimators, respectively, and have the following forms

$$\widehat{\text{Var}}(\widehat{\beta}_{\text{OLS}}) = \left(D' P_{X^\perp} D\right)^{-1} \widehat{\Sigma}_{11}, \quad \widehat{\text{Var}}(\widehat{\beta}_{\text{TSLS}}) = \left(D'(P_W - P_X)D\right)^{-1} \widehat{\Sigma}_{11}. \qquad (8)$$

The $\widehat{\Sigma}_{11}$ in equation (8) is the estimate of $\Sigma_{11}$ and can either be based on the OLS estimate, i.e. $\widehat{\Sigma}_{11} = \|Y - D\widehat{\beta}_{\text{OLS}} - X\widehat{\phi}_{\text{OLS}}\|_2^2/n$, or the TSLS estimate, i.e. $\widehat{\Sigma}_{11} = \|Y - D\widehat{\beta}_{\text{TSLS}} - X\widehat{\phi}_{\text{TSLS}}\|_2^2/n$.[5] Under $H_0$, both OLS and TSLS estimators of the variance $\Sigma_{11}$ are consistent. Also, under $H_0$, both OLS and TSLS estimators are consistent estimators of $\beta$, but the OLS estimator is more efficient than the TSLS estimator.

The asymptotic null distribution of the DWH test in equation (7) is a Chi-squared distribution with one degree of freedom ( the DWH test has an exact Chi-squared null distribution with one degree of freedom if $\Sigma_{11}$ is known). Hence for any $0 < \alpha < 1$, an asymptotically (or exactly if $\Sigma_{11}$ is known) level $\alpha$ test is given by

$$\text{Reject } H_0 \text{ if } \quad Q_{\text{DWH}} \geq \chi^2_\alpha(1).$$

Also, under the local alternative hypothesis,

$$H_0 : \Sigma_{12} = 0 \quad \text{versus} \quad H_2 : \Sigma_{12} = \frac{\Delta_1}{\sqrt{n}} \qquad (9)$$

for some constant $\Delta_1 \neq 0$, the asymptotic power of the DWH test is

$$\lim_{n \to \infty} \mathbf{P}(Q_{\text{DWH}} \geq \chi^2_\alpha(1)) = G\left(\alpha, \frac{\Delta_1 \sqrt{C(\mathcal{I})}}{\sqrt{\left(C(\mathcal{I}) + \frac{1}{p_z}\right) \Sigma_{11} \Sigma_{22}}}\right), \qquad (10)$$

where $G(\alpha, \cdot)$ is defined in equation (1); see Theorem 3 of the Supplementary Materials for a proof of equation (10). For textbook discussions on the DWH test, see Section 7.9 of Davidson and MacKinnon (1993) and Section 6.3.1 of Wooldridge (2010).

---

[5]The OLS and TSLS estimates of $\phi$ can be obtained as follows: $\widehat{\phi}_{\text{OLS}} = (X' P_{D^\perp} X)^{-1} X' P_{D^\perp} Y$ and $\widehat{\phi}_{\text{TSLS}} = (X' P_{\hat{D}^\perp} X)^{-1} X' P_{\hat{D}^\perp} Y$ where $\hat{D} = P_W D$.

## 3    The DWH Test with Many Covariates

We now consider the behavior of the DWH test in the presence of many covariates and/or instruments. Formally, suppose the number of covariates and instruments are growing with sample size $n$, that is, $p_x = p_x(n)$, $p_z = p_z(n)$ and $\lim_{n \to \infty} \min\{p_z, p_x\} = \infty$, so that $p = p_x + p_z$ and $n - p$ are increasing with respect to $n$. For this section only, we focus on the case where $p < n$ since the DWH test with OLS and TSLS estimators cannot be implemented when the sample size is smaller than the dimension of the model parameters; later sections, specifically Section 4, will consider endogeneity testing including both $p < n$ and $p \geq n$ settings. We assume a known $\Sigma_{11}$ for a cleaner technical exposition and to highlight the deficiencies of the DWH test that are not specific to estimating $\Sigma_{11}$, but specific to the form of the DWH test, the quadratic differencing of estimators in equation (7). However, the known $\Sigma_{11}$ assumption can be replaced by a consistent estimate of $\Sigma_{11}$. Theorem 1 characterizes the asymptotic behavior of the DWH test under this setting.

**Theorem 1.** *Suppose we have models (2) and (3) where $\Sigma_{11}$ is known, $W_{i\cdot}$ is a zero-mean multivariate Gaussian, the errors $\delta_i$ and $\epsilon_i$ are independent of $W_{i\cdot}$ and they are assumed to be bivariate normal. If $\sqrt{C(\mathcal{I})} \gg \sqrt{\log(n-p_x)/(n-p_x)p_z}$, for each $\alpha$, $0 < \alpha < 1$, the asymptotic Type I error of the DWH test under $H_0$ is controlled at $\alpha$, that is,*

$$\limsup_{n \to \infty} \mathbf{P}_\omega \left( |Q_{\text{DWH}}| \geq z_{\alpha/2} \right) = \alpha \ \text{ for any } \omega \text{ with corresponding } \Sigma_{12} = 0.$$

*Furthermore, for any $\omega$ with $\Sigma_{12} = \Delta_1/\sqrt{n}$, the asymptotic power of $Q_{\text{DWH}}$ satisfies*

$$\lim_{n \to \infty} \left| \mathbf{P}_\omega \left( Q_{\text{DWH}} \geq \chi^2_\alpha(1) \right) - G\left( \alpha, \frac{C(\mathcal{I}) \Delta_1 \sqrt{1 - \frac{p}{n}}}{\sqrt{\left( C(\mathcal{I}) + \frac{1}{n-p_x} \right) \left( C(\mathcal{I}) + \frac{1}{p_z} \right) \Sigma_{11} \Sigma_{22}}} \right) \right| = 0. \qquad (11)$$

Note that the convergence in Theorem 1 is pointwise convergence instead of uniform convergence. Theorem 1 states that the Type I error of the DWH test is actually controlled at the desired level $\alpha$ if one were to naively use it in the presence of many covariates and/or instruments and $\Sigma_{11}$ is known a priori. However, the power of the DWH test under the local alternative $H_2$ in equation (11) behaves differently in high dimensions than in low dimensions, as specified in equation (10). For example, if covariates and/or instruments are growing at $p/n \to 0$, equation (11) reduces to the usual power of the DWH test under low dimensional settings in equation (10). On the other hand, if covariates and/or instruments are growing at $p/n \to 1$, then the usual DWH test essentially has no power against any local alternative in $H_2$ since $G(\alpha, \cdot)$ in equation (11) equals $\alpha$ for any value of $\Delta_1$.

This phenomenon suggests that in the "middle ground" where $p/n \to c$, $0 < c < 1$, the usual DWH test will likely suffer in terms of power. As a concrete example, if $p_\text{x} = n/2$ and $p_\text{z} = n/3$ so that $p/n = 5/6$, then $G(\alpha, \cdot)$ in equation (11) reduces to

$$G\left(\alpha, \frac{C(\mathcal{I})\Delta_1}{\sqrt{2\left(C(\mathcal{I}) + \frac{2}{n}\right)\left(C(\mathcal{I}) + \frac{1}{p_z}\right)\Sigma_{11}\Sigma_{22}}}\right) \approx G\left(\alpha, \frac{1}{\sqrt{6}} \cdot \frac{\sqrt{C(\mathcal{I})}\Delta_1}{\sqrt{\left(C(\mathcal{I}) + \frac{1}{p_z}\right)\Sigma_{11}\Sigma_{22}}}\right)$$

where the approximation sign is for $n$ sufficiently large so that $C(\mathcal{I}) + 2/n \approx C(\mathcal{I})$. In this setting, the power of the DWH test is smaller than the power of the DWH test in equation (10) for the low dimensional setting. Section 5 also shows this phenomenon numerically.

Also, Theorem 1 provides some important guidelines for empiricists using the DWH test. First, Theorem 1 suggests that with modern cross-sectional data where the number of covariates may be very large, the DWH test should not be used to test endogeneity. Not only is the DWH test potentially incapable of detecting the presence of endogeneity under this scenario, but also an empiricist may be misled into an non-IV type of analysis, say the OLS or the Lasso, based on the result of the DWH test (Wooldridge, 2010). If the empiricist used a more powerful endogeneity test under this setting, he or she would have correctly concluded that there is endogeneity and used an IV analysis. Second, as discussed in Section 1, if empirical works add many covariates to make an IV more plausibly valid, one pays a price in terms of the power of the specification test; consequently, additional samples may be needed to get the desired level of power for detecting endogeneity.

Finally, we make two remarks about the regularity conditions in Theorem 1. First, Theorem 1 controls the growth of the concentration parameter $C(\mathcal{I})$ to be faster than $\log(n - p_{\text{x}})/(n - p_{\text{x}})p_{\text{z}}$. This growth condition is satisfied under the many instrument asymptotics of Bekker (1994) and the many weak instrument asymptotics of Chao and Swanson (2005) where $C(\mathcal{I})$ converges to a constant as $p_{\text{z}}/n \to c$ for some constant $c$. The weak instrument asymptotics of Staiger and Stock (1997) are not directly applicable to our growth condition on $C(\mathcal{I})$ because its asymptotics keeps $p_{\text{z}}$ and $p_{\text{x}}$ fixed. Second, we can replace the condition that $W_{i\cdot}$ is a zero-mean multivariate Gaussian in Theorem 1 by another condition used in high dimensional IV regression, for instance page 486 of Chernozhukov et al. (2015) where (i) the vector of instruments $Z_{i\cdot}$ is a linear model of $X_{i\cdot}$, i.e. $Z'_{i\cdot} = X'_{i\cdot} B + \bar{Z}'_{i\cdot}$, (ii) $\bar{Z}_{i\cdot}$ is independent of $X_{i\cdot}$, and (iii) $\bar{Z}_{i\cdot}$ is a multivariate normal distribution and the results in Theorem 1 will hold.

## 4  An Improved Endogeneity Test

Given that the DWH test for endogeneity may have low power in high dimensional settings, we present a simple and improved endogeneity test that has better power to detect endogeneity. In particular, our endogeneity test takes any popular estimator that is "well-behaved" for estimating reduced-form parameters (see Definition 2 for details) and applies a simple hard thresholding procedure to choose the most relevant instruments. We also stress that our endogeneity test is the first test capable of testing endogeneity if the number of parameters exceeds the sample size.

### 4.1  Well-Behaved Estimators

Consider the following reduced-form models

$$Y_i = Z'_{i\cdot}\Gamma + X'_{i\cdot}\Psi + \xi_i, \tag{12}$$

$$D_i = Z'_{i\cdot}\gamma + X'_{i\cdot}\psi + \epsilon_i. \tag{13}$$

The terms $\Gamma = \beta\gamma$ and $\Psi = \phi + \beta\psi$ are the parameters for the reduced-form model (12) and $\xi_i = \beta\epsilon_i + \delta_i$ is the reduced-form error term. The errors in the reduced-form models have the property that $\mathbf{E}(\xi_i|Z_{i.}, X_{i.}) = 0$ and $\mathbf{E}(\epsilon_i|Z_{i.}, X_{i.}) = 0$. Also, the covariance matrix of these error terms, denoted as $\Theta$, have the following forms: $\Theta_{11} = \text{Var}(\xi_i|Z_{i.}, X_{i.}) = \Sigma_{11} + 2\beta\Sigma_{12} + \beta^2\Sigma_{22}$, $\Theta_{22} = \text{Var}(\epsilon_i|Z_{i.}, X_{i.})$, and $\Theta_{12} = \text{Cov}(\xi_i, \epsilon_i|Z_{i.}, X_{i.}) = \Sigma_{12} + \beta\Sigma_{22}$.

As mentioned before, our improved endogeneity test does not require a specific estimator for the reduced-form parameters. Rather, any estimator that is well-behaved, as defined below, will be sufficient.

**Definition 2.** *Consider estimators $(\widehat{\gamma}, \widehat{\Gamma}, \widehat{\Theta}_{11}, \widehat{\Theta}_{22}, \widehat{\Theta}_{12})$ of the reduced-form parameters, $(\gamma, \Gamma, \Theta_{11}, \Theta_{22}, \Theta_{12})$ respectively, in equations (12) and (13). The estimators $(\widehat{\gamma}, \widehat{\Gamma}, \widehat{\Theta}_{11}, \widehat{\Theta}_{22}, \widehat{\Theta}_{12})$ are well-behaved estimators if they satisfy the two criterions below.*

*(W1) The reduced-form estimators of the coefficients $\widehat{\gamma}$ and $\widehat{\Gamma}$ satisfy*

$$\sqrt{n}\|(\widehat{\gamma} - \gamma) - \frac{1}{n}\widehat{V}'\epsilon\|_\infty = O_p\left(\frac{s\log p}{\sqrt{n}}\right), \quad \sqrt{n}\|(\widehat{\Gamma} - \Gamma) - \frac{1}{n}\widehat{V}'\xi\|_\infty = O_p\left(\frac{s\log p}{\sqrt{n}}\right). \tag{14}$$

*for some matrix $\widehat{V} = (\widehat{V}_{\cdot 1}, \cdots, \widehat{V}_{\cdot p_z})$ which is only a function of $W$ and satisfies*

$$\liminf_{n\to\infty} \mathbf{P}\left(c \le \min_{1\le j \le p_z} \frac{\|\widehat{V}_{\cdot j}\|_2}{\sqrt{n}} \le \max_{1\le j \le p_z} \frac{\|\widehat{V}_{\cdot j}\|_2}{\sqrt{n}} \le C, \ c\|\gamma\|_2 \le \frac{1}{\sqrt{n}}\|\sum_{j\in\mathcal{S}} \gamma_j \widehat{V}_{\cdot j}\|_2\right) = 1 \tag{15}$$

*for some constants $c > 0$ and $C > 0$.*

*(W2) The reduced-form estimators of the error variances, $\widehat{\Theta}_{11}$, $\widehat{\Theta}_{22}$, and $\widehat{\Theta}_{12}$ satisfy*

$$\sqrt{n}\max\left\{\left|\widehat{\Theta}_{11} - \frac{1}{n}\xi'\xi\right|, \left|\widehat{\Theta}_{12} - \frac{1}{n}\epsilon'\xi\right|, \left|\widehat{\Theta}_{22} - \frac{1}{n}\epsilon'\epsilon\right|\right\} = O_p\left(\frac{s\log p}{\sqrt{n}}\right). \tag{16}$$

There are many estimators of the reduced-form parameters in the literature that are well-behaved. Some examples of well-behaved estimators are listed below.

1. (OLS): In settings where $p$ is fixed or $p$ is growing with $n$ at a rate $p/n \to 0$, the OLS

estimators of the reduced-form parameters, i.e.

$$(\widehat{\Gamma}, \widehat{\Psi})' = (W'W)^{-1}W'Y \quad , (\widehat{\gamma}, \widehat{\psi})' = (W'W)^{-1}W'D,$$

$$\widehat{\Theta}_{11} = \frac{\left\|Y - Z\widehat{\Gamma} - X\widehat{\Psi}\right\|_2^2}{n-1} \quad , \widehat{\Theta}_{22} = \frac{\left\|D - Z\widehat{\gamma} - X\widehat{\psi}\right\|_2^2}{n-1}$$

$$\widehat{\Theta}_{12} = \frac{\left(Y - Z\widehat{\Gamma} - X\widehat{\Psi}\right)'\left(D - Z\widehat{\gamma} - X\widehat{\psi}\right)}{n-1}$$

trivially satisfy conditions for well-behaved estimators. Specifically, let $\widehat{V}' = (\frac{1}{n}W'W)^{-1}_{\mathcal{I}\cdot}W$. Then equation (14) holds because $(\widehat{\gamma} - \gamma) - \widehat{V}'\epsilon = 0$ and $(\widehat{\Gamma} - \Gamma) - \widehat{V}'\xi = 0$. Also, equation (15) holds because, in probability, $n^{-1/2}\|\widehat{V}_{\cdot j}\|_2 \to \Lambda_{jj}^{-1}$ and $n^{-1}\widehat{V}'\widehat{V} \to \Lambda_{\mathcal{I}\mathcal{I}}^{-1}$, thus satisfying (W1). Also, (W2) holds because $\|\widehat{\Gamma} - \Gamma\|_2^2 + \|\widehat{\Psi} - \Psi\|_2^2 = O_p(n^{-1})$ and $\|\widehat{\gamma} - \gamma\|_2^2 + \|\widehat{\psi} - \psi\|_2^2 = O_p(n^{-1})$, which implies equation (16) is going to zero at $n^{-1/2}$ rate.

2. (Debiased Lasso Estimators) In high dimensional settings where $p$ is growing with $n$ and often exceeds $n$, one of the most popular estimators for regression model parameters is the Lasso (Tibshirani, 1996). Unfortunately, the Lasso estimator and many penalized estimators do not satisfy the definition of a well-behaved estimator, specifically (W1), because penalized estimators are typically biased. Fortunately, recent works by Javanmard and Montanari (2014); van de Geer et al. (2014); Zhang and Zhang (2014) and Cai and Guo (2016) remedied this bias problem by doing a bias correction on the original penalized estimates.

More concretely, suppose we use the square root Lasso estimator by Belloni et al. (2011b),

$$\{\widetilde{\Gamma}, \widetilde{\Psi}\} = \underset{\Gamma \in \mathbb{R}^{p_z}, \Psi \in \mathbb{R}^{p_x}}{\operatorname{argmin}} \frac{\|Y - Z\Gamma - X\Psi\|_2}{\sqrt{n}} + \frac{\lambda_0}{\sqrt{n}}\left(\sum_{j=1}^{p_z}\|Z_{\cdot j}\|_2|\Gamma_j| + \sum_{j=1}^{p_x}\|X_{\cdot j}\|_2|\Psi_j|\right) \quad (17)$$

for the reduced-form model in equation (12) and

$$\{\widetilde{\gamma}, \widetilde{\psi}\} = \underset{\Gamma \in \mathbb{R}^{p_z}, \Psi \in \mathbb{R}^{p_x}}{\operatorname{argmin}} \frac{\|D - Z\gamma - X\psi\|_2}{\sqrt{n}} + \frac{\lambda_0}{\sqrt{n}} \left( \sum_{j=1}^{p_z} \|Z_{\cdot j}\|_2 |\gamma_j| + \sum_{j=1}^{p_x} \|X_{\cdot j}\|_2 |\psi_j| \right) \quad (18)$$

for the reduced-form model in equation (13). The term $\lambda_0$ in both estimation problems (17) and (18) represents the penalty term in the square root Lasso estimator and typically, the penalty is set at $\lambda_0 = \sqrt{a_0 \log p / n}$ for some constant $a_0$ slightly greater than 2, say 2.01 or 2.05. To transform the above penalized estimators in equations (17) and (18) into well-behaved estimators, we follow Javanmard and Montanari (2014) to debias the penalized estimators. Specifically, we solve $p_z$ optimization problems where the solution to each $p_z$ optimization problem, denoted as $\widehat{u}^{[j]} \in \mathbb{R}^p$, $j = 1, \ldots, p_z$, is

$$\widehat{u}^{[j]} = \underset{u \in \mathbb{R}^p}{\operatorname{argmin}} \frac{1}{n} \|Wu\|_2^2 \quad \text{s.t.} \quad \|\frac{1}{n} W'Wu - I_{\cdot j}\|_\infty \leq \lambda_n.$$

Typically, the tuning parameter $\lambda_n$ is chosen to be $12 M_1^2 \sqrt{\log p / n}$ where $M_1$ is defined as the largest eigenvalue of $\Lambda$. Define $\widehat{V}_{\cdot j} = W\widehat{u}^{[j]}$ and $\widehat{V} = (\widehat{V}_{\cdot 1}, \cdots, \widehat{V}_{\cdot p_z})$. Then, we can transform the penalized estimators in (17) and (18) into debiased, well-behaved estimators, $\widehat{\Gamma}$ and $\widehat{\gamma}$,

$$\widehat{\Gamma} = \widetilde{\Gamma} + \frac{1}{n}\widehat{V}'\left(Y - Z\widetilde{\Gamma} - X\widetilde{\Psi}\right), \quad \widehat{\gamma} = \widetilde{\gamma} + \frac{1}{n}\widehat{V}'\left(D - Z\widetilde{\gamma} - X\widetilde{\psi}\right). \quad (19)$$

Guo et al. (2016) showed that $\widehat{\Gamma}$ and $\widehat{\gamma}$ satisfy (W1). As for the error variances, following Belloni et al. (2011b), Sun and Zhang (2012) and Ren et al. (2013), we estimate the covariance terms $\Theta_{11}, \Theta_{22}, \Theta_{12}$ by

$$\widehat{\Theta}_{11} = \frac{\left\|Y - Z\widetilde{\Gamma} - X\widetilde{\Psi}\right\|_2^2}{n}, \widehat{\Theta}_{22} = \frac{\left\|D - Z\widetilde{\gamma} - X\widetilde{\psi}\right\|_2^2}{n}$$
$$\widehat{\Theta}_{12} = \frac{\left(Y - Z\widetilde{\Gamma} - X\widetilde{\Psi}\right)'\left(D - Z\widetilde{\gamma} - X\widetilde{\psi}\right)}{n}. \quad (20)$$

Lemma 3 of Guo et al. (2016) showed that the above estimators of $\widehat{\Theta}_{11}, \widehat{\Theta}_{22}$ and $\widehat{\Theta}_{12}$ in equation (20) satisfy (W2). In summary, the debiased Lasso estimators in equation (19) and the variance estimators in equation (20) are well-behaved estimators.

3. (One-Step and Orthogonal Estimating Equations Estimators) Recently, Chernozhukov et al. (2015) proposed the one-step estimator of the reduced-form coefficients, i.e.

$$\widehat{\Gamma} = \widetilde{\Gamma} + \frac{1}{n}\widetilde{\Lambda^{-1}}_{\mathcal{I},.}W^\intercal \left(Y - Z\widetilde{\Gamma} - X\widetilde{\Psi}\right), \widehat{\gamma} = \widetilde{\gamma} + \frac{1}{n}\widetilde{\Lambda^{-1}}_{\mathcal{I},.}W^\intercal \left(D - Z\widetilde{\gamma} - X\widetilde{\psi}\right).$$

where $\widetilde{\Gamma}, \widetilde{\gamma}$, and $\widetilde{\Lambda^{-1}}$ are initial estimators of $\Gamma$, $\gamma$ and $\Lambda^{-1}$, respectively. The initial estimators must satisfy conditions (18) and (20) of Chernozhukov et al. (2015) and many popular estimators like the Lasso or the square root Lasso satisfy these two conditions. Then, the arguments in Theorem 2.1 of van de Geer et al. (2014) showed that the one-step estimator of Chernozhukov et al. (2015) satisfies (W1). Relatedly, Chernozhukov et al. (2015) proposed estimators for the reduced-form coefficients based on orthogonal estimating equations and in Proposition 4 of Chernozhukov et al. (2015), the authors showed that the orthogonal estimating equations estimator is asymptotically equivalent to their one-step estimator.

For variance estimation, one can use the variance estimator in Belloni et al. (2011b), which reduces to the estimators in equation (20) and thus, satisfies (W2).

In short, the first part of our endogeneity test requires any estimator that is well-behaved and, as illustrated above, many estimators, such as the OLS in low dimensions and bias corrected penalized estimators in high dimensions, satisfy the criteria for a well-behaved estimator.

## 4.2 Estimating Relevant Instruments via Hard Thresholding

Once we have well-behaved estimators $(\widehat{\gamma}, \widehat{\Gamma}, \widehat{\Theta}_{11}, \widehat{\Theta}_{22}, \widehat{\Theta}_{12})$ satisfying Definition 2, the next step in our endogeneity test is finding IVs that are relevant, that is the set $\mathcal{S}$ in Definition 1 comprised of $\gamma_j \neq 0$. We do this by hard thresholding the estimate $\widehat{\gamma}$ by the dimension and the noise of $\widehat{\gamma}$.

$$\widehat{\mathcal{S}} = \left\{ j : |\widehat{\gamma}_j| \geq \frac{\sqrt{\widehat{\Theta}_{22}}\|\widehat{V}_{\cdot j}\|_2}{\sqrt{n}} \sqrt{\frac{a_0 \log \max\{p_z, n\}}{n}} \right\}. \tag{21}$$

The set $\widehat{\mathcal{S}}$ is an estimate of $\mathcal{S}$ and $a_0$ is some constant greater than 2; from our experience and like many Lasso problems, $a_0 = 2.01$ or $a_0 = 2.05$ works well in practice. The threshold in (21) is based on the noise level of $\widehat{\gamma}_j$ in equation (14) (represented by the term $n^{-1}\sqrt{\widehat{\Theta}_{22}}\|\widehat{V}_{\cdot j}\|_2$), adjusted by the dimensionality of the instrument size (represented by the term $\sqrt{a_0 \log \max\{p_z, n\}}$).

Using the estimated set $\widehat{\mathcal{S}}$ of relevant IVs leads to the estimates of $\Sigma_{12}$, $\Sigma_{11}$, and $\beta$,

$$\widehat{\Sigma}_{12} = \widehat{\Theta}_{12} - \widehat{\beta}\widehat{\Theta}_{22}, \quad \widehat{\Sigma}_{11} = \widehat{\Theta}_{11} + \widehat{\beta}^2\widehat{\Theta}_{22} - 2\widehat{\beta}\widehat{\Theta}_{12}, \quad \widehat{\beta} = \frac{\sum_{j \in \widehat{\mathcal{S}}} \widehat{\gamma}_j \widehat{\Gamma}_j}{\sum_{j \in \widehat{\mathcal{S}}} \widehat{\gamma}_j^2}. \tag{22}$$

Equation (22) provides us with the ingredients to construct our new test for endogeneity, which we denote as $Q$

$$Q = \frac{\sqrt{n}\widehat{\Sigma}_{12}}{\sqrt{\widehat{\mathrm{Var}}(\widehat{\Sigma}_{12})}} \quad \text{and} \quad \widehat{\mathrm{Var}}(\widehat{\Sigma}_{12}) = \widehat{\Theta}_{22}^2 \widehat{\mathrm{Var}}_1 + \widehat{\mathrm{Var}}_2 \tag{23}$$

where $\widehat{\mathrm{Var}}_1 = \widehat{\Sigma}_{11} \left\|\sum_{j \in \widehat{\mathcal{S}}} \widehat{\gamma}_j \widehat{V}_{\cdot j}/\sqrt{n}\right\|_2^2 / \left(\sum_{j \in \widehat{\mathcal{S}}} \widehat{\gamma}_j^2\right)^2$ and $\widehat{\mathrm{Var}}_2 = \widehat{\Theta}_{11}\widehat{\Theta}_{22} + \widehat{\Theta}_{12}^2 + 2\widehat{\beta}^2\widehat{\Theta}_{22}^2 - 4\widehat{\beta}\widehat{\Theta}_{12}\widehat{\Theta}_{22}$. Here, $\mathrm{Var}_1$ is the variance associated with estimating $\beta$ and $\mathrm{Var}_2$ is the variance associated with estimating $\Theta$.

A major difference between the original DWH test in equation (7) and our endogeneity test in equation (23) is that our endogeneity test directly estimates and tests the endogeneity parameter $\Sigma_{12}$ while the original DWH test implicitly tests for the endogeneity parameter by checking the quadratic distance between the OLS and TSLS estimators under the null hypothesis. More importantly, our endogeneity test efficiently uses the sparsity of the regression vectors while the DWH test does not incorporate such information. As shown in Section 4.3, our endogeneity test in this form where we make use of the sparsity information to estimate $\Sigma_{12}$ will have superior power in high dimension compared to the DWH test.

## 4.3 Properties of the New Endogeneity Test

We study the properties of our new test in high dimensional settings where $p$ is a function of $n$ and is allowed to be larger than $n$; note that this is a generalization of the setting discussed in Section 3 where $p < n$ because the DWH test is not feasible when $p \geq n$. Theorem 1 showed that the DWH test, while it controls Type I error at the desired level, may have low power, especially when the ratio of $p/n$ is close to 1. Theorem 2 shows that our new test $Q$ remedies this deficiency of the DWH test by having proper Type I error control and exhibiting better power than the DWH test.

**Theorem 2.** *Suppose we have models (2) and (3) where the errors $\delta_i$ and $\epsilon_i$ are independent of $W_i.$ and are assumed to be bivariate normal and we use a well-behaved estimator in our test statistic $Q$. If $\sqrt{C(\mathcal{S})} \gg s_{z1} \log p / \sqrt{n|\mathcal{V}|}$, and $\sqrt{s_{z1}} s \log p / \sqrt{n} \to 0$, then for any $\alpha$, $0 < \alpha < 1$, the asymptotic Type I error of $Q$ under $H_0$ is controlled at $\alpha$, that is,*

$$\lim_{n \to \infty} \mathbf{P}_w\left(|Q| \geq z_{\alpha/2}\right) = \alpha, \quad \text{for any } \omega \text{ with corresponding } \Sigma_{12} = 0. \tag{24}$$

*For any $\omega$ with $\Sigma_{12} = \Delta_1/\sqrt{n}$, the asymptotic power of $Q$ is*

$$\lim_{n \to \infty} \left| \mathbf{P}_\omega\left(|Q| \geq z_{\alpha/2}\right) - \mathbf{E}\left(G\left(\alpha, \frac{\Delta_1}{\sqrt{\Theta_{22}^2 \mathrm{Var}_1 + \mathrm{Var}_2}}\right)\right) \right| = 0, \tag{25}$$

*where $\mathrm{Var}_1 = \Sigma_{11} \left\| \sum_{j \in \mathcal{S}} \gamma_j \widehat{V}_{\cdot j}/\sqrt{n} \right\|_2^2 / \left(\sum_{j \in \mathcal{S}} \gamma_j^2\right)^2$ and $\mathrm{Var}_2 = \Theta_{11}\Theta_{22} + \Theta_{12}^2 + 2\beta^2 \Theta_{22}^2 - 4\beta \Theta_{12}\Theta_{22}$.*

In contrast to equation (11) that described the power of the usual DWH test in high dimensions, the term $\sqrt{1 - p/n}$ is absent in the power of our new endogeneity test $Q$ in equation (25). Specifically, under the local alternative $H_2$, our power is only affected by $\Delta_1$ while the power of the DWH test is affected by $\Delta_1 \sqrt{1 - p/n}$. Consequently, the power of our test $Q$ do not suffer from the growing dimensionality of $p$. For example, in the extreme case when $p/n \to 1$ and $C(\mathcal{S})$ is a constant, the power of the usual DWH test will be $\alpha$ while the power of our test $Q$ will always be greater than $\alpha$. For further validation, Section 5 numerically illustrates the discrepancies between the power of the two tests. Finally, we stress that in the case $p > n$, our test still has proper size and non-trivial power while the DWH test is not feasible in this setting.

With respect to the regularity conditions in Theorem 2, like Theorem 1, Theorem 2 controls the growth of the concentration parameter $C(\mathcal{S})$ to be faster than $s_{z1} \log p/\sqrt{n|\mathcal{S}|}$, with a minor discrepancy in the growth rate due to the differences between the set of relevant IVs, $\mathcal{S}$, and the set of candidate IVs, $\mathcal{I}$. But, similar to Theorem 1, this growth condition is satisfied under the many instrument asymptotics of Bekker (1994) and the many weak instrument asymptotics of Chao and Swanson (2005). Also, note that unlike the negative result in Theorem 1, the "positive" result in Theorem 2 is more general in that we do not require $W$ to be Gaussian and require $\Sigma_{11}$ to be known a priori. Instead, we only need the conditions of well-behaved estimators to hold. Also, we follow other high-dimensional inference works Javanmard and Montanari (2014); van de Geer et al. (2014); Zhang and Zhang (2014) in assuming independence and normality assumptions on the error terms $\delta_i$ and $\epsilon_i$, where such assumptions are made out of technicalities in establishing the distribution of test statistics in high dimensions. Finally, we remark that the expectation inside equation (25) is respect to $W$ and $\widehat{V}$ is a function of $W$.

## 4.4 An Extension: Endogeneity Test in High Dimensions with Possibly Invalid IVs

As discussed in Section 1, one of the motivations for having high dimensional covariates in empirical IV work is to avoid invalid instruments. While adding more covariates can potentially make instruments more plausibly valid, as demonstrated in Section 3, there is a price to pay with respect to the power of the DWH test. More importantly, even after conditioning on many covariates, some IVs may still be invalid and subsequent analysis, including the DWH test, assuming that all the IVs are valid after conditioning, can be seriously misleading. Inspired by these concerns, there has been a recent literature in estimation and inference of structural parameters in IV regression when invalid instruments are present(Guo et al., 2016; Kang et al., 2016; Kolesár et al., 2015). Our new endogeneity test $Q$ can be extended to handle the case of invalid instruments through the voting method proposed in Guo et al. (2016). The methodological and theoretical details are presented in Section 3.3 of the Supplementary Materials. To summarize the results, the extension of $Q$ to handle invalid instruments still controls the Type I error rate and has non-negligible power under high dimension with possibly invalid instruments.

# 5 Simulation and Data Example

## 5.1 Setup

We conduct a simulation study to investigate the performance of our new endogeneity test and the DWH test in high dimensional settings. Specifically, we generate data from models (2) and (3) in Section 2.2 with $n = 200$ or $300$, $p_z = 100$ and $p_x = 150$. The vector $W_{i\cdot}$ is a multivariate normal with mean zero and covariance $\Lambda_{ij} = 0.5^{|i-j|}$ for $1 \leq i,j \leq p$. We set the parameters as follows $\beta = 1$, $\phi = (0.6, 0.7, 0.8, \cdots, 1.5, 0, 0, \cdots, 0) \in \mathbb{R}^{p_x}$ so that $s_{x1} = 10$, and $\psi = (1.1, 1.2, 1.3, \cdots, 2.0, 0, 0, \cdots, 0) \in \mathbb{R}^{p_x}$ so that $s_{x2} = 10$. The relevant instruments are $\mathcal{S} = \{1, \ldots, 7\}$. Variance of the error terms are set to $\text{Var}(\delta_i) = \text{Var}(\epsilon_i) = 1.5$.

The parameters we vary in the simulation study are: the endogeneity level via $\text{Cov}(\delta_i, \epsilon_i)$, and IV strength via $\gamma$. For the endogeneity level, we set $\text{Cov}(\delta_i, \epsilon_i) = 1.5\rho$, where $\rho$ is varied and captures the level of endogeneity; a larger value of $|\rho|$ indicates a stronger correlation between the endogenous variable $D_i$ and the error term $\delta_i$. For IV strength, we set $\gamma_{\mathcal{S}} = K(1, 1, 1, 1, 1, 1, \rho_1)$ and $\gamma_{\mathcal{S}^C} = 0$, where $K$ is varied as a function of the concentration parameter (see below) and $\rho_1$ is either 0 or 0.2. Specifically, the value $K$ controls the global strength of instruments, with higher $|K|$ indicating strong instruments in a global sense. In contrast, the value $\rho_1$ controls the relative individual strength of instruments, specifically between the first six instruments in $\mathcal{S}$ and the seventh instrument. For example, $\rho_1 = 0.2$ implies that the seventh instrument's individual strength is only 20% of the first six instruments. Note that varying $\rho_1$ essentially stress-tests the thresholding step in our endogeneity test to numerically verify whether our testing procedure can handle relevant IVs with very small magnitudes of $\gamma$.

We specify $K$ as follows. Suppose we have a set of simulation parameters $\mathcal{S}, \rho_1, \Lambda$ and $\Sigma_{22}$. For each value of $100 \cdot C(\mathcal{S})$, we find the corresponding $K$ that satisfies $100 \cdot C(\mathcal{S}) = 100 \cdot K^2 \|\Lambda_{\mathcal{S}|\mathcal{S}^C}^{1/2}(1,1,1,1,1,1,\rho_1)\|_2^2/(7 \cdot 1.5)$. We vary $100 \cdot C(\mathcal{S})$ from 25 to 100, specifying $K$ for each value of $100 \cdot C(\mathcal{S})$.

For each simulation setting, we repeat the data generation 1000 times. For each simulation setting, we compare the power of our testing procedure $Q$ to the DWH test and the oracle DWH test where an oracle knows the support of the parameter vectors $\phi, \psi$ and $\gamma$. We set the desired $\alpha$ level for all three tests to be $\alpha = 0.05$.

## 5.2 Results

Table 1 and Figure 3 consider the high dimensional setting with $n = 200, 300$, $p_\mathrm{x} = 150$, and $p_\mathrm{z} = 100$. Table 1 measures the Type I error rate across three methods; for $n = 200$, the regular DWH test was not used since both the OLS and TSLS estimators are infeasible in this regime. We see a few clear trends in Table 1. First, generally speaking, all three methods control their Type I error around the desired $\alpha = 0.05$. Our proposed test has a slight upward bias of Type I error in some high dimensional settings with weak IV, i.e. where the $C$ value is around 25. But, the worst case upward bias is no more than 0.03 off from the target 0.05 and is within simulation error as $C$ gets larger. Additionally, as Figure 3 shows, the slight bias in Type I error in small $C$ regimes is offset by substantial power gains compared to the regular DWH test. Second, as the instrument gets stronger, both individually via $\rho_2$ and overall via $C$, the Type I error control generally gets better across all three methods, which is not surprising given the literature on strong instruments.

Figures 2 and 3 consider the power of our test $Q$, the regular DWH test, and the oracle DWH test in the high dimensional setting with $n = 200, 300, p_\mathrm{x} = 150$, and $p_\mathrm{z} = 100$. As predicted by Theorem 1, the regular DWH test suffers from low power, especially if the degree of endogeneity is around 0.25 where the gap between the regular DWH test and the oracle DWH test is the greatest across most simulation settings. In fact, even if the global strength of the IV increases, the DWH test still has low power. In contrast, as predicted from Theorem 2, our test $Q$ can handle $n \approx p$ or $n < p$. It also has uniformly better power than the regular DWH test across all degrees of endogeneity and across all simulation settings in the plot. Our test also achieves near-oracle performance as the global instrument strength grows.

|   |     | Weak    |       |        | Strong  |       |        |
|---|-----|---------|-------|--------|---------|-------|--------|
| C | n   | Regular | Ours  | Oracle | Regular | Ours  | Oracle |
| 25  | 300 | 0.040 | 0.079 | 0.034 | 0.061 | 0.048 | 0.038 |
|     | 200 | NA    | 0.080 | 0.054 | NA    | 0.075 | 0.054 |
| 50  | 300 | 0.049 | 0.046 | 0.032 | 0.043 | 0.065 | 0.048 |
|     | 200 | NA    | 0.072 | 0.055 | NA    | 0.069 | 0.050 |
| 75  | 300 | 0.053 | 0.059 | 0.044 | 0.043 | 0.062 | 0.048 |
|     | 200 | NA    | 0.065 | 0.038 | NA    | 0.063 | 0.048 |
| 100 | 300 | 0.067 | 0.055 | 0.048 | 0.050 | 0.064 | 0.044 |
|     | 200 | NA    | 0.057 | 0.045 | NA    | 0.049 | 0.045 |

Table 1: Empirical Type I error when $p_x = 150$ and $p_z = 100$ after 1000 simulations. The value $n$ represents the sample size and $\alpha = 0.05$. "Regular," "Ours," and "Oracle" represent the regular DWH test, the proposed test $(Q)$, and the oracle DWH test, respectively. "Weak", and "Strong" represent the cases when $\rho_1 = 0.2$ and $\rho_1 = 0$, respectively. $C$ represents the overall strength of the instruments, as measured by $100 \cdot C(\mathcal{S})$. NA indicates not applicable.

In summary, all the simulation results indicate that our endogeneity test controls Type I error and is a much better alternative to the regular DWH test in high dimensional settings, with near-optimal performance with respect to the oracle. Our test is also capable of handling the regime $n < p$. In the supplementary materials, we also conduct low dimensional simulations and show that all three tests, the oracle DWH test, the regular DWH test, and our proposed test behave identically with respect to power and Type I error control.

## 5.3 Data Example

To highlight the usefulness of the proposed test statistic $Q$, specifically its ability to run DWH test in dimensions where $n < p$, we re-analyze a high dimensional data analysis done in Belloni et al. (2012, 2014). Specifically, the outcome $Y$ is the log of average Case-Shiller home price index and the endogenous variable $D$ is the number of federal appellate court decisions that were against seizure of property via eminent domain. There are $n = 183$ individuals and $p_z = 147$ instruments which are derived from indicators that represent the random assignment of judges to different cases, characteristics of judges, and other interactions. Additionally, there are $p_x = 71$ exogenous variables that describe the type of cases, number of court decisions, circuit specific and time-specific effects; see Belloni et al. (2012) and Belloni et al. (2014) for more details about the instruments and the exogenous variables. We use the code provided in Belloni et al. (2012) to replicate the data set.

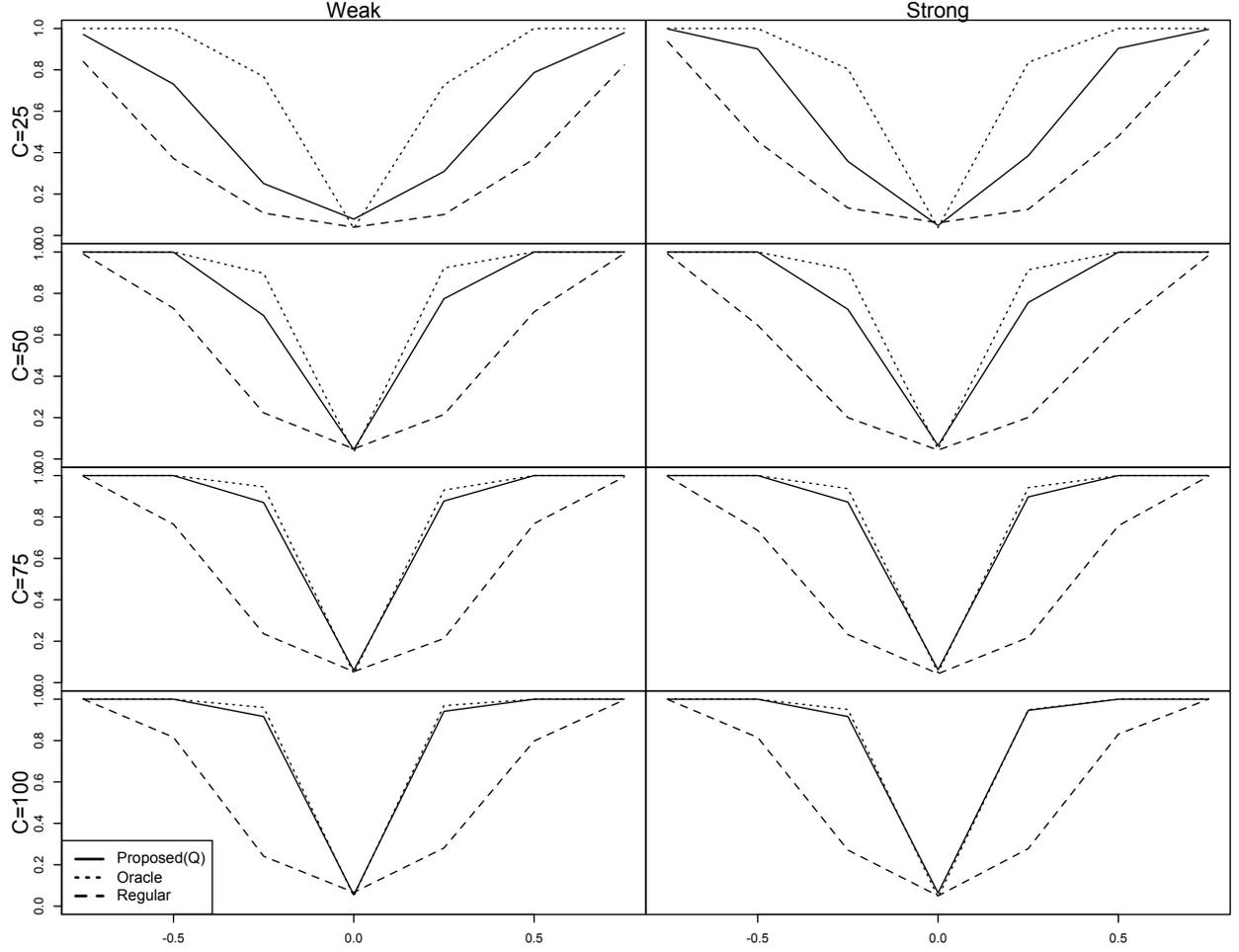

Figure 1: Power of endogeneity tests when $n = 300$, $p_x = 150$ and $p_z = 100$. The $x$-axis represents the endogeneity $\rho$ and the y-axis represents the empirical power over 1000 simulations. Each line represents a particular test's empirical power over various values of the endogeneity, where the solid line, the dashed line and the dotted line represent the proposed test ($Q$), the regular DWH test and the oracle DWH test, respectively. The columns represent the individual IV strengths, with column names "Weak" and "Strong" denoting the cases when $\rho_1 = 0.2$, and $\rho_1 = 0$, respectively. The rows represent the overalls strength of the instruments, as measured by $100 \cdot C(\mathcal{S})$.

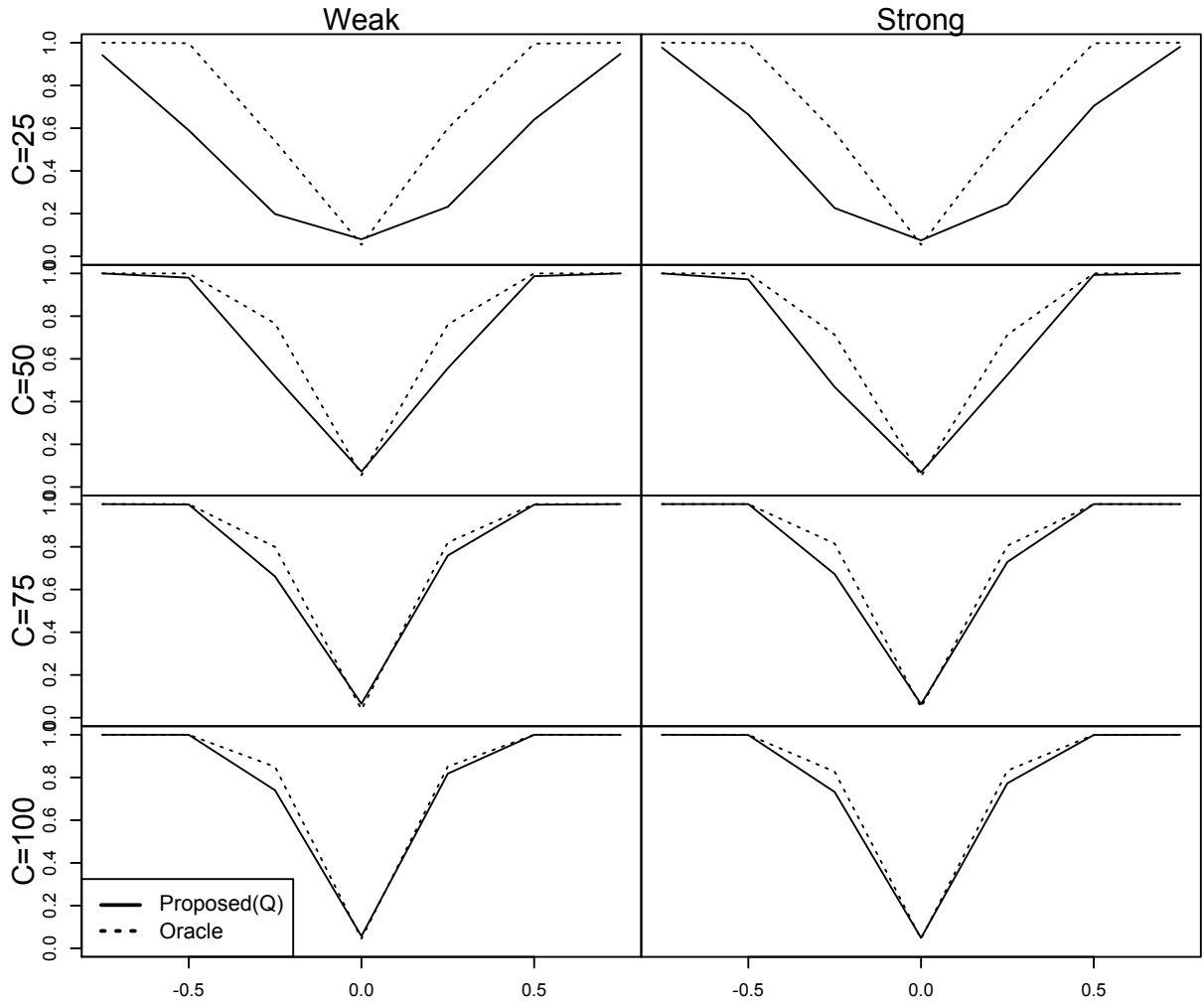

Figure 2: Power of endogeneity tests when $n = 200$, $p_x = 150$ and $p_z = 100$. The $x$-axis represents the endogeneity $\rho$ and the y-axis represents the empirical power over 1000 simulations. Each line represents a particular test's empirical power over various values of the endogeneity, where the solid line and the dotted line represent the proposed test ($Q$) and the oracle DWH test, respectively. The columns represent the individual IV strengths, with column names "Weak" and "Strong" denoting the cases when $\rho_1 = 0.2$, and $\rho_1 = 0$, respectively. The rows represent the overalls strength of the instruments, as measured by $100 \cdot C(\mathcal{S})$.

Because $n < p$, the DWH test or other tests for endogeneity cannot be used. Consequently, investigators are forced to remove covariates and/or instruments to run their usual specification test. For example, in our analysis, we drop the covariates and use the AER package (Kleiber and Zeileis, 2008), which is a popular R package to run IV analysis, to run the DWH test. The package reports back that the p-value for the DWH test is 0.683.

In contrast, our new test $Q$ allows data where $n < p$. As such, we are not forced to remove covariates from the original analysis when we run our test $Q$ on this data. Our test reports the p-value for the Q test is 0.21, meaning that there is not evidence for the number of federal appellate court decisions against seizure of property or eminent domain being endogenous. Unlike the DWH test, our test was able to accommodate these high dimensional covariates rather than dropping them from the analysis.

## 6 Conclusion and Discussion

In this paper, we showed that the popular DWH test, while being able to control Type I error, can have low power in high dimensional settings. We propose a simple and improved endogeneity test to remedy the low power of the DWH test by modifying popular reduced-form parameters with a thresholding step. We also show that this modification leads to drastically better power than the DWH test in high dimensional settings.

For empirical work, the results in the paper suggest that one should be cautious in interpreting high $p$-values produced by the DWH test in IV regression settings when many covariates and/or instruments are present. In particular, as shown in Section 3, in modern data settings with a potentially large number of covariates and/or instruments, the DWH test may declare that there is no endogeneity in the structural model, even if endogeneity is truly present. Our proposed test, which is a simple modification of the popular estimators for reduced-forms parameters, does not suffer from this problem, as it achieves near-oracle performance to detect endogeneity, and can even handle general settings when $n < p$ and invalid IVs are present.

## Acknowledgments

The research of Hyunseung Kang was supported in part by NSF Grant DMS-1502437. The research of T. Tony Cai was supported in part by NSF Grants DMS-1403708 and DMS-1712735, and NIH Grant R01 GM-123056. The research of Dylan S. Small was supported in part by NSF Grant SES-1260782.

# Supplement to "Testing Endogeneity with High Dimensional Covariates"

Zijian Guo[1], Hyunseung Kang[2], T. Tony Cai[3], and Dylan S. Small[3]

[1]Department of Statistics and Biostatistics, Rutgers University
[2]Department of Statistics, University of Wisconsin-Madison
[3]Department of Statistics, The Wharton School, University of Pennsylvania

**Abstract**

This note summarizes the supplementary materials to the paper "Testing Endogeneity with High Dimensional Covariates". In Section 1, we present extended simulation studies for the low dimensional setting. In Section 2, we show that the DWH test fails in the presence of Invalid IVs. In Section 3, we discuss both method and theory for endogeneity test in high dimension with invalid IVs. In Section 4, we present technical proofs for Theorems 1, 2, 3, 4 and 5 and the proofs of technical lemmas.

## 1 Simulation for Low Dimensions

For the low dimensional case, we generate data from models the same models as the high dimensional simulations, except we have $p_{\rm z} = 9$ instruments, $p_{\rm x} = 5$ covariates, and $n = 1000$ samples. The parameters of the models are: $\beta = 1$, $\phi = (0.6, 0.7, 0.8, 0.9, 1.0) \in \mathbb{R}^5$ and $\psi = (1.1, 1.2, 1.3, 1.4, 1.5) \in \mathbb{R}^5$. We see that the three comparators, the regular DWH test, the oracle DWH test, and our test are very similar with respect to power and Type I error control.

*Address for correspondence: Zijian Guo, Department of Statistics and Biostatistics, Rutgers University, USA. Phone: (848)445-2690. Fax: (732)445-3428. Email: zijguo@stat.rutgers.edu.

## 2 Failure of the DWH Test in the Presence of Invalid IVs

While the DWH test performs as expected when all the instruments are valid, in practice, some instruments may be invalid and consequently, the DWH test can be a highly misleading assessment of the hypotheses (6). In Theorem 3, we show that the Type I error of the DWH test can be greater than the nominal level for a wide range of IV configurations in which some IVs are invalid; we assume a known $\Sigma_{11}$ in Theorem 3 for a cleaner technical exposition and to highlight the impact that invalid IVs have on the size and power of the DWH test, but the known $\Sigma_{11}$ can be replaced by a consistent estimate of $\Sigma_{11}$. We also show that the power of the DWH test under the local alternative $H_2$ in equation (9) can be shifted.

**Theorem 3.** *Suppose we have models (2) and (3) with a known $\Sigma_{11}$. If $\pi = \Delta_2/n^k$ where $\Delta_2$ is a fixed constant and $0 \leq k < \infty$, then for any $\alpha$, $0 < \alpha < 1$, we have the following asymptotic phase-transition behaviors of the DWH test for different values of $k$.*

a. *$0 \leq k < 1/2$: The asymptotic Type I error of the DWH test under $H_0$ is 1, i.e.*

$$\omega \in H_0 : \lim_{n \to \infty} \mathbf{P}\left(Q_{\text{DWH}} \geq \chi^2_\alpha(1)\right) = 1 \tag{26}$$

*and the asymptotic power of the DWH test under $H_2$ is 1.*

b. *$k = 1/2$: The asymptotic Type I error of the DWH test under $H_0$ is*

$$\omega \in H_0 : \lim_{n \to \infty} \mathbf{P}\left(Q_{\text{DWH}} \geq \chi^2_\alpha(1)\right) = G\left(\alpha, \frac{\frac{1}{p_z}\gamma'\Lambda_{\mathcal{I}|\mathcal{I}^c}\Delta_2}{\sqrt{C(\mathcal{I})\left(C(\mathcal{I}) + \frac{1}{p_z}\right)\Sigma_{11}\Sigma_{22}}}\right) \geq \alpha, \tag{27}$$

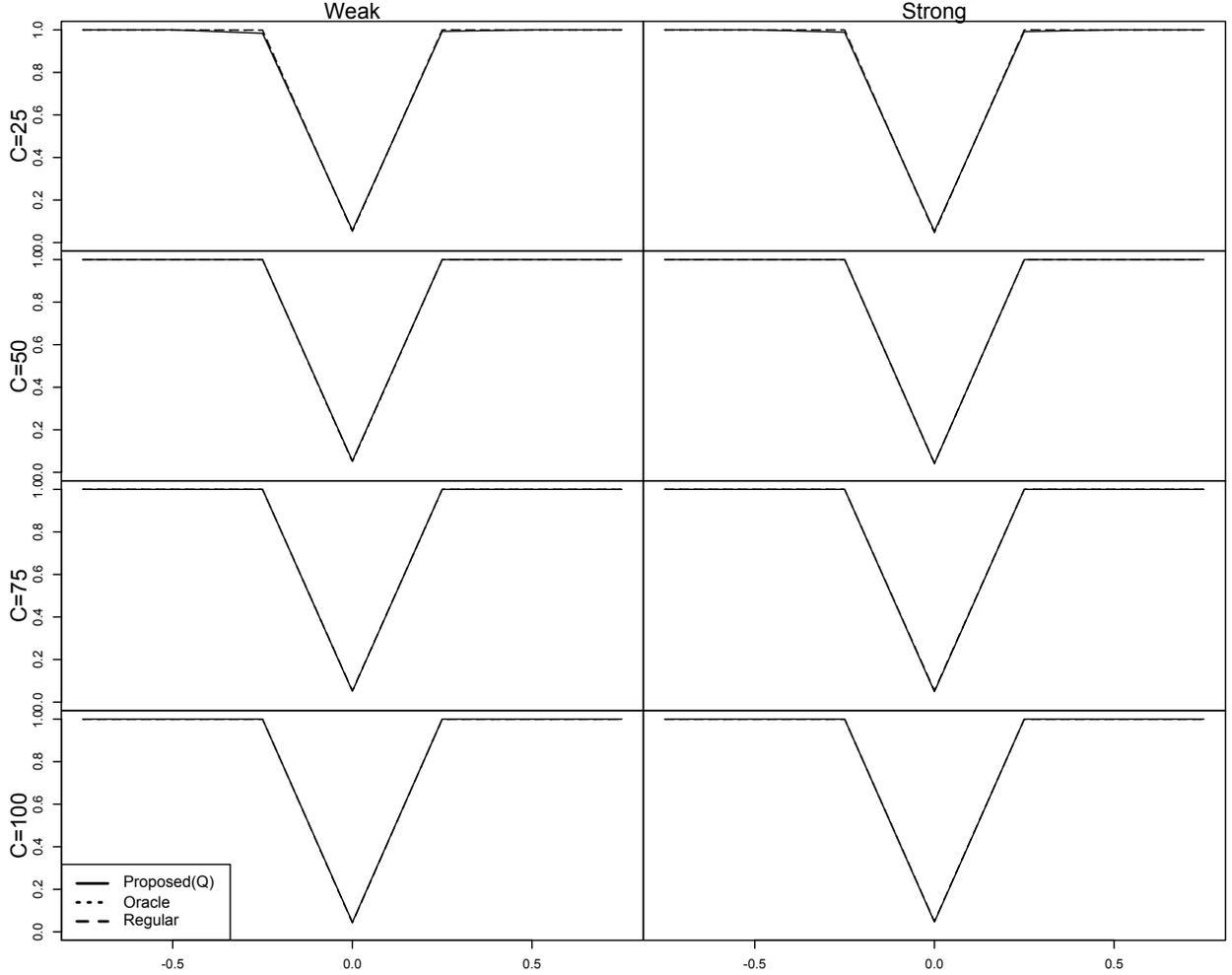

Figure 3: Power of endogeneity tests when $n = 1000$, $p_{\text{x}} = 5$ and $p_{\text{z}} = 9$. The $x$-axis represents the endogeneity $\rho$ and the y-axis represents the empirical power over 1000 simulations. Each line represents a particular test's empirical power over various values of the endogeneity, where the solid line, the dashed line and the dotted line represent the proposed test $(Q)$, the regular DWH test and the oracle DWH test, respectively. The columns represent the individual IV strengths, with column names "Weak" and "Strong" denoting the cases when $\rho_1 = 0.2$, and $\rho_1 = 0$, respectively. The rows represent the overalls strength of the instruments, as measured by $100 \cdot C(\mathcal{S})$.

*and the asymptotic power of the DWH test under $H_2$ is*

$$\omega \in H_2 : \lim_{n \to \infty} \mathbf{P}\left(Q_{\text{DWH}} \geq \chi^2_\alpha(1)\right)$$

$$= G\left(\alpha, \frac{\frac{1}{p_z}\gamma' \Lambda_{\mathcal{I}|\mathcal{I}^c}\Delta_2}{\sqrt{C(\mathcal{I})\left(C(\mathcal{I}) + \frac{1}{p_z}\right)\Sigma_{11}\Sigma_{22}}} + \frac{\Delta_1\sqrt{C(\mathcal{I})}}{\sqrt{\left(C(\mathcal{I}) + \frac{1}{p_z}\right)\Sigma_{11}\Sigma_{22}}}\right), \quad (28)$$

*where $G(\alpha, \cdot)$ is defined in* (1).

c. $1/2 < k < \infty$: *The asymptotic Type I error of the DWH test is $\alpha$, i.e.*

$$\omega \in H_0 : \lim_{n \to \infty} \mathbf{P}\left(Q_{\text{DWH}} \geq \chi^2_\alpha(1)\right) = \alpha \quad (29)$$

*and the asymptotic power of the DWH test under $H_2$ is equivalent to equation* (10).

Theorem 3 presents the asymptotic behavior of the DWH test under a wide range of settings for the invalid IVs as represented by $\pi$. For example, when the instruments are invalid in the sense that their deviation from valid IVs (i.e. $\pi = 0$) to invalid IVs (i.e. $\pi \neq 0$) is at rates slower than $n^{-1/2}$, say $\pi = \Delta_2 n^{-1/4}$ or $\pi = \Delta_2$, equation (26) states that the DWH will always have Type I error and power that reach 1. In other words, if some IVs, or even a single IV, are moderately (or strongly) invalid in the sense that they have moderate (or strong) direct effects on the outcome above the usual noise level of the model error terms at $n^{-1/2}$, then the DWH test will always reject the null hypothesis of no endogeneity even if there is truly no endogeneity present; essentially, the DWH test behaves equivalently to a test that never looks at the data and always rejects the null.

Next, suppose the instruments are invalid in the sense that their deviation from valid IVs to invalid IVs are exactly at $n^{-1/2}$ rate, also referred to as the Pitman drift.[1] This is the phase-transition point of the DWH test's Type I error as the error moves from 1 in equation (26) to $\alpha$ in equation (29). Under this type of invalidity, equation (27) shows that the Type I error of the DWH test depends on some factors, most prominently the factor $\gamma' \Lambda_{\mathcal{I}|\mathcal{I}^c} \Delta_2$. The factor $\gamma' \Lambda_{\mathcal{I}|\mathcal{I}^c} \Delta_2$ has been discussed in the literature, most recently by Kolesár et al. (2015) within the context of invalid IVs. Specifically, Kolesár et al. (2015) studied the case where $\Delta_2 \neq 0$ so that there are invalid IVs, but $\gamma' \Lambda_{\mathcal{I}|\mathcal{I}^c} \Delta_2 = 0$, which essentially amounted to saying that the IVs' effect on the endogenous variable $D$ via $\gamma$ is orthogonal to their direct effects on the outcome via $\Delta_2$; see Assumption 5 of Section 3 in Kolesár et al. (2015) for details. Under their scenario, if $\gamma' \Lambda_{\mathcal{I}|\mathcal{I}^c} \Delta_2 = 0$, then the DWH test will have the desired size $\alpha$. However, if $\gamma' \Lambda_{\mathcal{I}|\mathcal{I}^c} \Delta_2$ is not exactly zero, which will most likely be the case in practice, then the Type I error of the DWH test will always be larger than $\alpha$ and we can compute the exact deviation from $\alpha$ by using equation (27). Also, equation (28) computes the power under $H_2$ in the $n^{-1/2}$ setting, which again depends on the magnitude and direction of $\gamma' \Lambda_{\mathcal{I}|\mathcal{I}^c} \Delta_2$. For example, if there is only one instrument and that instrument has average negative effects on both $D$ and $Y$, the overall effect on the power curve will be a positive shift away from the case of valid IVs (i.e. $\pi = 0$). Regardless, under the $n^{-1/2}$ invalid IV regime, the DWH test will always have size that is at least as large as $\alpha$ if invalid IVs are present.

---

[1] Fisher (1967) and Newey (1985) used this type of $n^{-1/2}$ asymptotic argument to study misspecified econometrics models, specifically Section 2, equation (2.3) of Fisher (1967) and Section 2, Assumption 2 of Newey (1985). More recently, Hahn and Hausman (2005) and Berkowitz et al. (2012) used the $n^{-1/2}$ asymptotic framework in their respective works to study plausibly exogenous variables.

Theorem 3 also shows that instruments' strength, as measured by the population concentration parameter $C(\mathcal{I})$ in equation (5), impacts the Type I error rate of the DWH test when the IVs are invalid at the $n^{-1/2}$ rate. Specifically, if $\pi = \Delta_2 n^{-1/2}$ and the instruments are strong so that the concentration parameter $C(\mathcal{I})$ is large, then the deviation from $\alpha$ will be relatively minor even if $\gamma' \Lambda_{\mathcal{I}|\mathcal{I}^c} \Delta_2 \neq 0$. This phenomenon has been mentioned in previous work, most notably Bound et al. (1995) and Angrist et al. (1996) where strong instruments can lessen the undesirable effects caused by invalid IVs.

Finally, if the instruments are invalid in the sense that their deviation from $\pi = 0$ is faster than $n^{-1/2}$, say $\pi = \Delta n^{-1}$, then equation (29) shows that the DWH test maintains its desired size. To put this invalid IV regime in context, if the instruments are invalid at $n^{-k}$ where $k > 1/2$, the convergence toward $\pi = 0$ is faster than the usual convergence rate of a sample mean from an i.i.d. sample towards a population mean. Also, this type of deviation is equivalent to saying that the invalid IVs are very weakly invalid and essentially act as if they are valid because the IVs are below the noise level of the model error terms at $n^{-1/2}$. Consequently, the DWH test is not impacted by these type of IVs with respect to size and power.

The overall implication of Theorem 3 is that whenever there is a concern for instrument validity, the results of the DWH test in practice should be scrutinized, especially when the DWH test produces low $p$-values. In particular, our theorem shows that the DWH test will only have correct size, (i) when the invalid IVs essentially behave as valid IVs asymptotically so that $\pi$'s rate toward zero is faster than usual mean convergence or (ii) when the IVs' effects on the endogenous variables are completely orthogonal to each other. In all other settings, the Type I error of the DWH test will often be larger than $\alpha$ and consequently, the DWH test will tend to over-reject the null more frequently than it should, even if a single invalid IV is present. In fact, the low $p$-value of the DWH test may mislead empiricists about the true presence of endogeneity; the endogeous variable may actually be exogenous and the low $p$-value may be entirely an artifact due to invalid IVs.

## 3 Endogeneity Test in high dimensions with Invalid IVs

### 3.1 Model

In this line of work[2] ,the invalid instruments are represented as direct effects between the instruments and the outcome in equation (2), i.e.

$$Y_i = D_i \beta + Z'_{i.} \pi + X'_{i.} \phi + \delta_i, \quad E(\delta_i \mid Z_{i.}, X_{i.}) = 0 \tag{30}$$

If $\pi = 0$ in model (30), the model (30) reduces to the usual instrumental variables regression model in equation (2) with one endogenous variable, $p_x$ exogenous covariates, and $p_z$ instruments, all of which are assumed to be valid. On the other hand, if $\pi \neq 0$ and the support of $\pi$ is unknown a priori, the instruments may have a direct effect on the outcome, thereby violating the exclusion restriction (Angrist et al., 1996; Imbens and Angrist, 1994), without knowing, a priori, which are invalid and valid (Conley et al., 2012; Kang et al., 2016; Murray, 2006). In short, the support of $\pi$ allows us to distinguish a valid instrument, i.e. $\pi_j = 0$ from an invalid one, i.e. $\pi_j \neq 0$.

---

[2]Works by Berkowitz et al. (2012); Fisher (1966, 1967); Guggenberger (2012); Hahn and Hausman (2005); Newey (1985) and Caner (2014) also considered properties of IV estimators or, more broadly, generalized method of moments estimators (GMM)s when there are local deviations from validity to invalidity.Andrews (1999) and Andrews and Lu (2001) considered selecting valid instruments within the context of GMMs. Small (2007) approached the invalid instrument problem via a sensitivity analysis. Conley et al. (2012) proposed various strategies, including union-bound correction, sensitivity analysis, and Bayesian analysis, to deal with invalid instruments. Liao (2013) and Cheng and Liao (2015) considered the setting where there is, a priori, a known set of valid instruments and another set of instruments that may not be valid.

## 3.2 Method

Despite the presence of invalid IVs, our new endogeneity test can handle this case by using an additional thresholding procedure outlined in Section 3.3 of Guo et al. (2016a) to estimate $\pi$ in the model (30). Specifically, we take each IV $j$ that are estimated to be relevant, i.e. $j \in \widehat{\mathcal{S}}$, and we define $\widehat{\beta}^{[j]}$ to be a "pilot" estimate of $\pi$ by using this IV and dividing the reduced-form parameter estimates, i.e. $\widetilde{\pi}^{[j]} = \widehat{\Gamma} - \widehat{\beta}^{[j]}\widehat{\gamma}$ where $\widehat{\beta}^{[j]} = \widehat{\Gamma}_j/\widehat{\gamma}_j$. We also define $\widehat{\pi}^{[j]}$ to be a pilot estimate of $\pi$ using this $j$th instrument's estimate of $\beta$, i.e. $\widetilde{\pi}^{[j]} = \widehat{\Gamma} - \widehat{\beta}^{[j]}\widehat{\gamma}$, and $\widehat{\Sigma}_{11}^{[j]}$ to be the pilot estimate of $\Sigma_{11}$, i.e. $\widehat{\Sigma}_{11}^{[j]} = \widehat{\Theta}_{11} + (\widehat{\beta}^{[j]})^2\widehat{\Theta}_{22} - 2\widehat{\beta}^{[j]}\widehat{\Theta}_{12}$. Then, for each $\widetilde{\pi}^{[j]}$ in $j \in \widehat{\mathcal{S}}$, we threshold each element of $\widetilde{\pi}^{[j]}$ to create the thresholded estimate $\widehat{\pi}^{[j]}$,

$$\widehat{\pi}_k^{[j]} = \widetilde{\pi}_k^{[j]} \mathbf{1}\left(k \in \widehat{\mathcal{S}} \ \cap \ |\widetilde{\pi}_k^{[j]}| \geq a_0 \sqrt{\widehat{\Sigma}_{11}^{[j]}} \frac{\|\widehat{V}_{\cdot k} - \frac{\widehat{\gamma}_k}{\widehat{\gamma}_j}\widehat{V}_{\cdot j}\|_2}{\sqrt{n}} \sqrt{\frac{\log \max(p_z, n)}{n}}\right) \quad (31)$$

for all $1 \leq k \leq p_z$. Each thresholded estimate $\widehat{\pi}^{[j]}$ is obtained by looking at the elements of the un-thresholded estimate, $\widetilde{\pi}^{[j]}$, and examining whether each element exceeds the noise threshold (represented by the term $n^{-1}\sqrt{\widehat{\Sigma}_{11}^{[j]}}\|\widehat{V}_{\cdot k} - \frac{\widehat{\gamma}_k}{\widehat{\gamma}_j}\widehat{V}_{\cdot j}\|_2$), adjusted for the multiplicity of the selection procedure (represented by the term $a_0\sqrt{\log \max(p_z, n)}$). Among the $|\widehat{\mathcal{S}}|$ candidate estimates of $\pi$ based on each relevant instrument in $\widehat{\mathcal{S}}$, i.e. $\widehat{\pi}^{[j]}$, we choose $\widehat{\pi}^{[j]}$ with the most valid instruments, i.e. we choose $j^* \in \widehat{\mathcal{S}}$ where $j^* = \operatorname{argmin} \|\widehat{\pi}^{[j]}\|_0$; if there is a non-unique solution, we choose $\widehat{\pi}^{[j]}$ with the smallest $\ell_1$ norm, the closest convex norm of $\ell_0$. Subsequently, we can estimate the set of valid and relevant IVs, denoted as $\widehat{\mathcal{V}} \subseteq \mathcal{I}$, as those elements of $\widehat{\pi}^{[j^*]}$ that are zero,

$$\widehat{\mathcal{V}} = \widehat{\mathcal{S}} \setminus \text{supp}\left(\widehat{\pi}^{[j^*]}\right). \tag{32}$$

and estimate $\beta$ as

$$\widehat{\beta} = \frac{\sum_{j \in \widehat{\mathcal{V}}} \widehat{\gamma}_j \widehat{\Gamma}_j}{\sum_{j \in \widehat{\mathcal{V}}} \widehat{\gamma}_j^2}. \tag{33}$$

The endogeneity test that is robust to invalid IVs has the same form as equation (23), except we use the set $\widehat{\mathcal{V}}$ instead of $\widehat{\mathcal{S}}$ and the estimate $\widehat{\beta}_E$ instead of $\widehat{\beta}$. We denote this endogeneity test as $Q_E$.

### 3.3 Properties of $Q_E$

We analyze the properties of the endogeneity test $Q_E$, which can handle invalid instruments as well as high dimensional instruments and covariates, even when $p > n$. Let $\mathcal{V} = \{j \in \mathcal{I} \mid \pi_j = 0, \gamma_j \neq 0\}$. We make the following assumptions that essentially control the behavior of selecting relevant and invalid IVs. We denote the assumption as "IN" since the assumption is specific to the case when invalid IVs are present.

(IN1) (50% Rule) The number of valid IVs is more than half of the number of non-redundant IVs, that is $|\mathcal{V}| > \frac{1}{2}|\mathcal{S}|$.

(IN2) (Individual IV Strength) Among IVs in $\mathcal{S}$, we have $\min_{j \in \mathcal{S}} |\gamma_j| \geq \delta_{\min} \gg \sqrt{\log p / n}$.

(IN3) (Strong violation) Among IVs in the set $\mathcal{S} \setminus \mathcal{V}$, we have

$$\min_{j \in \mathcal{S} \setminus \mathcal{V}} \left|\frac{\pi_j}{\gamma_j}\right| \geq \frac{12(1+|\beta|)}{\delta_{\min}} \sqrt{\frac{M_1 \log \max\{p_z, n\}}{n \lambda_{\min}(\Theta)}}. \tag{34}$$

In a nutshell, Assumption (IN1) states that if the number of invalid instruments is not too large, then we can use the observed data to separate the invalid IVs from valid IVs, without knowing a priori which IVs are valid or invalid. Assumption (IN1) is a relaxation of the assumption typical in IV settings where all the IVs are assumed to be valid a priori so that $|\mathcal{V}| = p_z$ and (IN1) holds automatically. In particular, Assumption (IN1) entertains the possibility that some IVs may be invalid, so $|\mathcal{V}| < p_z$, but without knowing a priori which IVs are invalid, i.e. the exact set $\mathcal{V}$. Assumption (IN1) is also the generalization of the 50% rule in Han (2008) and Kang et al. (2016) in the presence of redundant IVs. Also, Kang et al. (2016) showed that this type of proportion-based assumption is crucial for identification of model parameters when instrument validity is uncertain.

Assumption (IN2) requires individual IV strength to be bounded away from zero. This assumption is needed to rule out IVs that are asymptotically weak. We also show in the simulation studies presented in the supplementary materials that (IN2) is largely unnecessary for our test to have proper size and have good power. Also, in the literature, (IN2) is similar to the "beta-min" condition assumption in high dimensional linear regression without IVs (Bühlmann and Van De Geer, 2011; Fan and Li, 2001; Wainwright, 2007; Zhao and Yu, 2006), with the exception that this condition is not imposed on our inferential quantity of interest, the endogeneity parameter $\Sigma_{12}$. Next, Assumption (IN3) requires the ratio $\pi_j/\gamma_j$ for invalid IVs to be large. This assumption is needed to correctly select valid IVs in the presence of possibly invalid IVs and this sentiment is echoed in the model selection literature by Leeb and Pötscher (2005) who pointed out that "in general no model selector can be uniformly consistent for the most parsimonious true model" and hence the post-model-selection inference is generally non-uniform (or uniform within a limited class of models). Specifically, for any IV with a small, but non-zero $|\pi_j/\gamma_j|$, such a weakly invalid IV is hard to distinguish from valid IVs where $\pi_j/\gamma_j = 0$. If a weakly invalid IV is mistakenly declared as valid, the bias from this mistake is of the order $\sqrt{\log p_z/n}$, which has consequences, not for consistency of the point estimation of $\Sigma_{12}$, but for a $\sqrt{n}$ inference of $\Sigma_{12}$.

If all the instruments are valid, like the setting described in the majority of this paper where the IVs are valid conditional on many covariates, we do not need Assumptions (IN1)-(IN3) to make any claims about the proposed endogeneity test. However, in the presence of potentially invalid IVs that can grow in dimension, assumptions (IN1)-(IN3) are needed to control the behavior of the invalid IVs asymptotically and to characterize the asymptotic behavior of $Q_E$.

**Theorem 4.** *Suppose we have models (2) and (3) where the errors $\delta_i$ and $\epsilon_i$ are independent of $W_i$. and are assumed to be bivariate normal but some instruments may be invalid, i.e. $\pi \neq 0$, and Assumptions* (IN1)-(IN3) *hold. If $\sqrt{C(\mathcal{V})} \gg s_{z1} \log p/\sqrt{n|\mathcal{V}|}$, and $\sqrt{s_{z1}} s \log p/\sqrt{n} \to 0$, then for any $\alpha$, $0 < \alpha < 1$, he asymptotic Type I error of $Q$ under $H_0$ is controlled at $\alpha$, that is,*

$$\lim_{n \to \infty} \mathbf{P}_w \left( |Q_E| \geq z_{\alpha/2} \right) = \alpha, \quad \text{for any } \omega \text{ with corresponding } \Sigma_{12} = 0. \tag{35}$$

*For any $\omega$ with $\Sigma_{12} = \Delta_1/\sqrt{n}$, the asymptotic power of $Q_E$ is*

$$\lim_{n \to \infty} \left| \mathbf{P}_\omega \left( |Q_E| \geq z_{\alpha/2} \right) - \mathbf{E} \left( G \left( \alpha, \frac{\Delta_1}{\sqrt{\Theta_{22}^2 \text{Var}_1 + \text{Var}_2}} \right) \right) \right| = 0, \tag{36}$$

*where* $\text{Var}_1 = \Sigma_{11} \left\| \sum_{j \in \mathcal{V}} \gamma_j \widehat{V}_{\cdot j}/\sqrt{n} \right\|_2^2 / \left( \sum_{j \in \mathcal{S}} \gamma_j^2 \right)^2$ *and* $\text{Var}_2 = \Theta_{11}\Theta_{22} + \Theta_{12}^2 + 2\beta^2 \Theta_{22}^2 - 4\beta\Theta_{12}\Theta_{22}$.

Theorem 4 shows that our new test $Q_E$ controls Type I error at the desired level $\alpha$. Also, Theorem 4 states that the power of $Q_E$ is similar to the power of $Q$ that knows exactly which instruments are valid and relevant. In short, our test $Q_E$ is adaptive to the knowledge about instrument validity and can achieve similar level of performance as the test $Q$ that knows about instrument validity a priori.

Finally, like Theorem 2, Theorem 4 controls the growth of the concentration parameter $C(\mathcal{V})$ to be faster than $s_{z1} \log p/\sqrt{n|\mathcal{V}|}$, with a minor discrepancy in the growth rate due to the differences between the sets $\mathcal{V}$ and $\mathcal{S}$. But, as before, this growth condition is satisfied under the many instrument asymptotics of Bekker (1994) and the many weak instrument asymptotics of Chao and Swanson (2005). Also, like Theorem 4, the regularity conditions on $s, s_{z1}, p, n$ are the same as those from Theorem 2.

# 4 Proof

## 4.1 Proof of Theorem 3

<u>Proof of (29)</u> By the assumption $(\delta_i, \epsilon_i) \sim N\left(0, \begin{pmatrix} \Sigma_{11} & \Sigma_{12} \\ \Sigma_{21} & \Sigma_{22} \end{pmatrix}\right)$, we have the following decomposition,

$$\delta_i = \frac{\Sigma_{12}}{\Sigma_{22}}\epsilon_i + \tau_i, \tag{37}$$

where $\tau_i$ is independent of $\epsilon_i$. By plugging (37) into (2) in the main paper, we have

$$Y_i = D_i\beta + Z'_{i.}\pi + X'_{i.}\phi + \frac{\Sigma_{12}}{\Sigma_{22}}\epsilon_i + \tau_i.$$

Let $\sigma_\tau^2$ denote the variance of $\tau_i$ and then $\sigma_\tau = \sqrt{\Sigma_{11} - \frac{\Sigma_{12}^2}{\Sigma_{22}}}$. Define

$$a_0(n) = \frac{\sigma_\tau}{\sqrt{\Sigma_{11}}} - 1 = \frac{-\frac{\Sigma_{12}^2}{\Sigma_{11}\Sigma_{22}}}{\sqrt{1 - \frac{\Sigma_{12}^2}{\Sigma_{11}\Sigma_{22}}} + 1}. \tag{38}$$

By the definition $\Sigma_{12} = \frac{\Delta}{\sqrt{n}}$, we have

$$|a_0(n)| \leq C\frac{1}{n}. \tag{39}$$

By the explicit expression of $\widehat{\beta}_{\text{OLS}}$ and $\widehat{\beta}_{\text{TSLS}}$,

$$\widehat{\beta}_{\text{OLS}} = \beta + (D'P_{X^\perp}D)^{-1}D'P_{X^\perp}(Z,\epsilon)\begin{pmatrix} \pi \\ \frac{\Sigma_{12}}{\Sigma_{22}} \end{pmatrix} + (D'P_{X^\perp}D)^{-1}D'P_{X^\perp}\tau$$

and

$$\widehat{\beta}_{\text{TSLS}} = \beta + (D'(P_W - P_X)D)^{-1}D'(P_W - P_X)(Z,\epsilon)\begin{pmatrix} \pi \\ \frac{\Sigma_{12}}{\Sigma_{22}} \end{pmatrix}$$

$$+ (D'(P_W - P_X)D)^{-1}D'(P_W - P_X)\tau.$$

we obtain the following decomposition of the difference $\widehat{\beta}_{\text{TSLS}} - \widehat{\beta}_{\text{OLS}}$,

$$\widehat{\beta}_{\text{TSLS}} - \widehat{\beta}_{\text{OLS}} = \left((D'(P_W - P_X)D)^{-1} - (D'P_{X^\perp}D)^{-1}\right) D'P_{X^\perp}Z\pi$$
$$+ \left((D'(P_W - P_X)D)^{-1}D'(P_W - P_X) - (D'P_{X^\perp}D)^{-1}D'P_{X^\perp}\right)\epsilon\frac{\Sigma_{12}}{\Sigma_{22}} \quad (40)$$
$$+ \left((D'(P_W - P_X)D)^{-1}D'(P_W - P_X) - (D'P_{X^\perp}D)^{-1}D'P_{X^\perp}\right)\tau.$$

In the following, we analyze the three terms in the above decomposition,

1. Conditioning on $\epsilon$ and $W$, we have

$$L_1 = \frac{\left((D'(P_W - P_X)D)^{-1}D'(P_W - P_X) - (D'P_{X^\perp}D)^{-1}D'P_{X^\perp}\right)\tau}{\sqrt{(D'(P_W - P_X)D)^{-1}\sigma_\tau^2 - (D'P_{X^\perp}D)^{-1}\sigma_\tau^2}} \sim N(0, 1). \quad (41)$$

2. By the assumption $\text{Cov}(W_i) = \Lambda$, $\text{Cov}\begin{pmatrix} \delta_i \\ \epsilon_i \end{pmatrix} = \Sigma$ and weak law of large number, we have

$$\frac{1}{n}Z'Z \xrightarrow{p} \Lambda_{zz}, \quad \frac{1}{n}X'Z \xrightarrow{p} \Lambda_{xz}, \quad \frac{1}{n}X'X \xrightarrow{p} \Lambda_{xx},$$
$$\frac{1}{n}\epsilon'Z \xrightarrow{p} 0, \quad \frac{1}{n}\epsilon'X \xrightarrow{p} 0, \quad \frac{1}{n}\epsilon'W \xrightarrow{p} 0, \quad \frac{1}{n}\epsilon'\epsilon \xrightarrow{p} \Sigma_{22}.$$

Hence, we have

$$(\frac{1}{n}D'P_{X^\perp}D)^{-1} \xrightarrow{p} \left(\gamma'\Lambda_{\mathcal{I}|\mathcal{I}^c}\gamma + \Sigma_{22}\right)^{-1}, \quad (\frac{1}{n}D'(P_W - P_X)D)^{-1} \xrightarrow{p} \left(\gamma'\Lambda_{\mathcal{I}|\mathcal{I}^c}\gamma\right)^{-1},$$
(42)

$$\frac{1}{n}D'P_{X^\perp}Z\pi \xrightarrow{p} \gamma'\Lambda_{\mathcal{I}|\mathcal{I}^c}\pi, \quad \frac{1}{n}D'(P_W - P_X)\epsilon \xrightarrow{p} 0 \quad \frac{1}{n}D'P_{X^\perp}\epsilon \xrightarrow{p} \Sigma_{22}, \quad (43)$$

By (42) and (43) and the parametrization $\Sigma_{12} = \frac{\Delta_1}{\sqrt{n}}$, we have

$$L_2 = \frac{\left((D'(P_W - P_X)D)^{-1}D'(P_W - P_X) - (D'P_{X^\perp}D)^{-1}D'P_{X^\perp}\right)\epsilon\frac{\Sigma_{12}}{\Sigma_{22}}}{\sqrt{(D'(P_W - P_X)D)^{-1}\sigma_\tau^2 - (D'P_{X^\perp}D)^{-1}\sigma_\tau^2}}$$
$$\xrightarrow{p} L_2^* = -\Delta_1\sqrt{\frac{\gamma'\Lambda_{\mathcal{I}|\mathcal{I}^c}\gamma}{\Sigma_{11}\Sigma_{22}\left(\gamma'\Lambda_{\mathcal{I}|\mathcal{I}^c}\gamma + \Sigma_{22}\right)}}.$$
(44)

3. By (42) and (43) and the parametrization $\pi = \frac{\Delta_2}{\sqrt{n}}$ where $\Delta_2$ is fixed vector, we have

$$L_3 = \frac{\left((D'(P_W - P_X)D)^{-1} - (D'P_{X^\perp}D)^{-1}\right) D'P_{X^\perp}Z\pi}{\sqrt{(D'(P_W - P_X)D)^{-1}\sigma_\tau^2 - (D'P_{X^\perp}D)^{-1}\sigma_\tau^2}}$$
$$\xrightarrow{p} L_3^* = \frac{\gamma' \Lambda_{\mathcal{I}|\mathcal{I}^c} \Delta_2 \sqrt{\Sigma_{22}}}{\sqrt{(\gamma' \Lambda_{\mathcal{I}|\mathcal{I}^c} \gamma + \Sigma_{22})(\gamma' \Lambda_{\mathcal{I}|\mathcal{I}^c} \gamma)} \sqrt{\Sigma_{11}}}. \tag{45}$$

Together with (40), we derive the general power curve as follows,

$$\mathbf{P}\left((L_1 + L_2 + L_3)^2 \geq \chi_\alpha^2(1)\right)$$
$$=\mathbf{P}\left(L_1 + L_2 + L_3 \geq \sqrt{\chi_\alpha^2(1)}\right) + \mathbf{P}\left(L_1 + L_2 + L_3 \leq -\sqrt{\chi_\alpha^2(1)}\right)$$
$$=\mathbf{P}\left(L_1 \geq \sqrt{\chi_\alpha^2(1)} - L_2 - L_3\right) + \mathbf{P}\left(L_1 \leq -\sqrt{\chi_\alpha^2(1)} - L_2 - L_3\right)$$
$$=\mathbf{E}_{W,\epsilon}\left(\mathbf{P}\left(L_1 \geq \sqrt{\chi_\alpha^2(1)} - L_2 - L_3 \mid W, \epsilon\right) + \mathbf{P}\left(L_1 \leq -\sqrt{\chi_\alpha^2(1)} - L_2 - L_3 \mid W, \epsilon\right)\right).$$

By (41), conditioning on $W$ and $\epsilon$, we have

$$\mathbf{P}\left(L_1 \geq \sqrt{\chi_\alpha^2(1)} - L_2 \mid W, \epsilon\right) = 1 - \Psi\left(\frac{\sqrt{\chi_\alpha^2(1)} - L_2 - L_3}{1 + a_0(n)}\right),$$
$$\mathbf{P}\left(L_1 \leq -\sqrt{\chi_\alpha^2(1)} - L_2 \mid W, \epsilon\right) = \Psi\left(\frac{-\sqrt{\chi_\alpha^2(1)} - L_2 - L_3}{1 + a_0(n)}\right).$$

Combined with (39), (41), (48) and (49), we establish (29). The type I error control (27) follows from (28) with taking $\Delta_2 = 0$.

<u>Proof of (26) and (29)</u> For the case $0 \leq k < \frac{1}{2}$, we apply the similar argument as that of (28) to establish (26) and the only difference is that

$$\frac{L_3}{\gamma' \Lambda_{\mathcal{I}|\mathcal{I}^c} \Delta_2} \xrightarrow{p} \frac{\sqrt{\Sigma_{22}}}{\sqrt{(\gamma' \Lambda_{\mathcal{I}|\mathcal{I}^c} \gamma + \Sigma_{22})(\gamma' \Lambda_{\mathcal{I}|\mathcal{I}^c} \gamma)} \sqrt{\Sigma_{11}}}. \tag{46}$$

As $\gamma' \Lambda_{\mathcal{I}|\mathcal{I}^c} \Delta_2 \to \infty$, we establish (26). For the case $k > \frac{1}{2}$, we apply the similar argument as that of (28) and we can establish (29) with the fact $L_3 \xrightarrow{p} 0$.

## 4.2 Proof of Theorem 1

By (40), we have the following expression of $\widehat{\beta}_{\text{TSLS}} - \widehat{\beta}_{\text{OLS}}$,

$$\widehat{\beta}_{\text{TSLS}} - \widehat{\beta}_{\text{OLS}} = \left((D'(P_W - P_X)D)^{-1}D'(P_W - P_X) - (D'P_{X^\perp}D)^{-1}D'P_{X^\perp}\right)\epsilon\frac{\Sigma_{12}}{\Sigma_{22}}$$
$$+ \left((D'(P_W - P_X)D)^{-1}D'(P_W - P_X) - (D'P_{X^\perp}D)^{-1}D'P_{X^\perp}\right)\tau.$$

Hence, the test statistic $Q_{\text{DWH}}$ has the following expression,

$$Q_{\text{DWH}} = \frac{\left(\widehat{\beta}_{\text{TSLS}} - \widehat{\beta}_{\text{OLS}}\right)^2}{(D'(P_W - P_X)D)^{-1}\Sigma_{11} - (D'P_{X^\perp}D)^{-1}\Sigma_{11}} = (L_1 + L_2)^2, \qquad (47)$$

where

$$L_1 = \frac{\left((D'(P_W - P_X)D)^{-1}D'(P_W - P_X) - (D'P_{X^\perp}D)^{-1}D'P_{X^\perp}\right)\tau}{\sqrt{(D'(P_W - P_X)D)^{-1}\sigma_\tau^2 - (D'P_{X^\perp}D)^{-1}\sigma_\tau^2}} \times \frac{\sigma_\tau}{\sqrt{\Sigma_{11}}} \qquad (48)$$

and

$$\begin{aligned}L_2 &= \frac{\left((D'(P_W - P_X)D)^{-1}D'(P_W - P_X) - (D'P_{X^\perp}D)^{-1}D'P_{X^\perp}\right)\epsilon\frac{\Sigma_{12}}{\Sigma_{22}}}{\sqrt{(D'(P_W - P_X)D)^{-1}\Sigma_{11} - (D'P_{X^\perp}D)^{-1}\Sigma_{11}}} \\ &= \frac{\left((D'(P_W - P_X)D)^{-1} - (D'P_{X^\perp}D)^{-1}\right)D'(P_W - P_X)\epsilon\frac{\Sigma_{12}}{\Sigma_{22}}}{\sqrt{(D'(P_W - P_X)D)^{-1}\Sigma_{11} - (D'P_{X^\perp}D)^{-1}\Sigma_{11}}} \\ &\quad - \frac{(D'P_{X^\perp}D)^{-1}D'\left(P_{X^\perp} - (P_W - P_X)\right)\epsilon\frac{\Sigma_{12}}{\Sigma_{22}}}{\sqrt{(D'(P_W - P_X)D)^{-1}\Sigma_{11} - (D'P_{X^\perp}D)^{-1}\Sigma_{11}}} \\ &= \frac{\sqrt{n}\Sigma_{12}}{\Sigma_{22}\sqrt{\Sigma_{11}}}\sqrt{(\frac{1}{n}D'(P_W - P_X)D)^{-1} - (\frac{1}{n}D'P_{X^\perp}D)^{-1}}\frac{1}{n}D'(P_W - P_X)\epsilon \\ &\quad - \frac{\sqrt{n}\Sigma_{12}}{\Sigma_{22}\sqrt{\Sigma_{11}}}\frac{(\frac{1}{n}D'P_{X^\perp}D)^{-1}\frac{1}{n}D'\left(P_{X^\perp} - (P_W - P_X)\right)\epsilon}{\sqrt{(\frac{1}{n}D'(P_W - P_X)D)^{-1} - (\frac{1}{n}D'P_{X^\perp}D)^{-1}}}.\end{aligned} \qquad (49)$$

In the following, we further decompose $L_2$. Since $W_{i\cdot}$ is a zero-mean multivariate Gaussian, we have $Z'_{i\cdot} = X'_{i\cdot}\left(\Lambda_{xx}^{-1}\Lambda_{xz}\right) + \bar{Z}'_{i\cdot}$ where $\bar{Z}_{i\cdot}$ is independent of $X_{i\cdot}$ and $\bar{Z}_{i\cdot}$ is of mean 0 and covariance matrix $\Lambda_{\mathcal{I}|\mathcal{I}^c}$. Hence, we have $D = \bar{Z}\gamma + X\left(\psi + \Lambda_{xx}^{-1}\Lambda_{xz}\gamma\right) + \epsilon$ and $\bar{Z}_{i\cdot}\gamma \sim N(0, \gamma'\Lambda_{\mathcal{I}|\mathcal{I}^c}\gamma)$. We further decompose $\frac{1}{n}D'(P_W - P_X)\epsilon$ as,

$$\begin{aligned}\frac{1}{n}D'(P_W - P_X)\epsilon &= \frac{1}{n}\left(\bar{Z}\gamma + X\left(\psi + \Lambda_{xx}^{-1}\Lambda_{xz}\gamma\right) + \epsilon\right)'\left(W\left(W'W\right)^{-1}W' - X\left(X'X\right)^{-1}X'\right)\epsilon \\ &= \frac{1}{n}\gamma'\bar{Z}'\left(I - X\left(X'X\right)^{-1}X'\right)\epsilon + \frac{1}{n}\epsilon'\left(W\left(W'W\right)^{-1}W' - X\left(X'X\right)^{-1}X'\right)\epsilon.\end{aligned}$$

Together with (49), we have

$$L_2 = L_{2,1} \times (L_{2,2} + L_{2,3}) - L_{2,4}. \tag{50}$$

where

$$L_{2,1} = \frac{\sqrt{n}\Sigma_{12}}{\Sigma_{22}\sqrt{\Sigma_{11}}}\sqrt{(\frac{1}{n}D'(P_W - P_X)D)^{-1} - (\frac{1}{n}D'P_{X^\perp}D)^{-1}}, \quad L_{2,2} = \frac{1}{n}\gamma'\bar{Z}'\left(I - X(X'X)^{-1}X'\right)\epsilon,$$

$$L_{2,3} = \frac{1}{n}\epsilon'\left(W(W'W)^{-1}W' - X(X'X)^{-1}X'\right)\epsilon, \quad L_{2,4} = \frac{\sqrt{n}\Sigma_{12}}{\Sigma_{22}\sqrt{\Sigma_{11}}}\frac{(\frac{1}{n}D'P_{X^\perp}D)^{-1}\frac{1}{n}D'(P_{X^\perp} - (P_W - P_X))\epsilon}{\sqrt{(\frac{1}{n}D'(P_W - P_X)D)^{-1} - (\frac{1}{n}D'P_{X^\perp}D)^{-1}}}.$$

In the following, we derive the power function. By (47), we derive the power function as follows,

$$\mathbf{P}\left((L_1 + L_2)^2 \geq \chi_\alpha^2(1)\right) = \mathbf{P}\left(L_1 \geq \sqrt{\chi_\alpha^2(1)} - L_2\right) + \mathbf{P}\left(L_1 \leq -\sqrt{\chi_\alpha^2(1)} - L_2\right)$$
$$= \mathbf{E}_{W,\epsilon}\left(\mathbf{P}\left(L_1 \geq \sqrt{\chi_\alpha^2(1)} - L_2 \mid W, \epsilon\right) + \mathbf{P}\left(L_1 \leq -\sqrt{\chi_\alpha^2(1)} - L_2 \mid W, \epsilon\right)\right). \tag{51}$$

By (49) and (48), conditioning on $W$ and $\epsilon$, we have

$$\mathbf{P}\left(L_1 \geq \sqrt{\chi_\alpha^2(1)} - L_2 \mid W, \epsilon\right) = 1 - \Psi\left(\frac{\sqrt{\chi_\alpha^2(1)} - L_2}{1 + a_0(n)}\right),$$
$$\mathbf{P}\left(L_1 \leq -\sqrt{\chi_\alpha^2(1)} - L_2 \mid W, \epsilon\right) = \Psi\left(\frac{-\sqrt{\chi_\alpha^2(1)} - L_2}{1 + a_0(n)}\right). \tag{52}$$

By (50), (51) and (52), we have

$$\mathbf{P}\left((L_1 + L_2)^2 \geq \chi_\alpha^2(1)\right) = \mathbf{E}_{W,\epsilon}\left(1 - \Psi\left(\frac{-\sqrt{\chi_\alpha^2(1)} - L_{2,1}(L_{2,2} + L_{2,3}) + L_{2,4}}{1 + a_0(n)}\right)\right)$$
$$+ \mathbf{E}_{W,\epsilon}\Psi\left(\frac{\sqrt{\chi_\alpha^2(1)} - L_{2,1}(L_{2,2} + L_{2,3}) + L_{2,4}}{1 + a_0(n)}\right). \tag{53}$$

In the following, we further approximate the power curve in (53). We first approximate the terms $L_{2,i}$ by $L_{2,i}^*$ for $i = 1, 3, 4$,

$$L_{2,1}^* = \frac{\sqrt{n}\Sigma_{12}}{\Sigma_{22}\sqrt{\Sigma_{11}}} \sqrt{\frac{\frac{n-p}{n}\Sigma_{22}}{\left(\frac{n-p_x}{n}\left(\gamma'\Lambda_{\mathcal{I}|\mathcal{I}^c}\gamma + \Sigma_{22}\right)\right)\left(\frac{n-p_x}{n}\left(\gamma'\Lambda_{\mathcal{I}|\mathcal{I}^c}\gamma\right) + \frac{p_z}{n}\Sigma_{22}\right)}}$$

$$= \frac{\sqrt{n}\Sigma_{12}}{\Sigma_{22}\sqrt{\Sigma_{11}}} \sqrt{\frac{\frac{n(n-p)}{(n-p_x)^2 p_z^2 \Sigma_{22}}}{\left(\frac{\gamma'\Lambda_{\mathcal{I}|\mathcal{I}^c}\gamma}{p_z \Sigma_{22}} + \frac{1}{n-p_x}\right)\left(\frac{\gamma'\Lambda_{\mathcal{I}|\mathcal{I}^c}\gamma}{p_z \Sigma_{22}} + \frac{1}{p_z}\right)}}$$

$$L_{2,3}^* = \frac{p_z}{n}\Sigma_{22}, \quad L_{2,4}^* = \frac{\sqrt{n}\Sigma_{12}}{\Sigma_{22}\sqrt{\Sigma_{11}}} \sqrt{\frac{(1-\frac{p}{n})\Sigma_{22}\left(\frac{\gamma'\Lambda_{\mathcal{I}|\mathcal{I}^c}\gamma}{p_z \Sigma_{22}} + \frac{1}{n-p_x}\right)}{\frac{\gamma'\Lambda_{\mathcal{I}|\mathcal{I}^c}\gamma}{p_z \Sigma_{22}} + \frac{1}{p_z}}}.$$

The following Lemma characterizes the difference between $L_{2,i}$ and $L_{2,i}^*$ for $i = 1, 3, 4$.

**Lemma 1.** *Define*

$$a_1(n) = \frac{L_{2,1}}{L_{2,1}^*} - 1, \quad \text{and} \quad a_2(n) = \frac{L_{2,3}}{L_{2,3}^*} - 1, \quad \text{and} \quad a_3(n) = \frac{L_{2,4}}{L_{2,4}^*} - 1. \tag{54}$$

*Then there exists an event $\mathcal{A}$ such that*

$$\mathbf{P}(\mathcal{A}) \geq 1 - (\min\{n - p, p_z\})^{-c}. \tag{55}$$

*On the event $\mathcal{A}$, there exists some positive constant $C$ such that*

$$\max\{|a_1(n)|, |a_3(n)|\} \leq C \left(\frac{\log p_z}{p_z} + \frac{\log(n-p)}{n-p} + \sqrt{\frac{\log(n-p_x)}{n-p_x}}\sqrt{\frac{\Sigma_{22}}{\gamma'\Lambda_{\mathcal{I}|\mathcal{I}^c}\gamma}}\right), \tag{56}$$

$$|a_2(n)| \leq C\frac{\log p_z}{p_z}. \tag{57}$$

In the following, we use Lemma 1 and calculate the approximation error to the exact power function (53),

$$
\left| \mathbf{E}_{W,\epsilon} \Psi \left( \frac{\sqrt{\chi_\alpha^2(1)} - L_{2,1}(L_{2,2} + L_{2,3}) + L_{2,4}}{1 + a_0(n)} \right) - \Psi \left( \sqrt{\chi_\alpha^2(1)} - L_{2,1}^* L_{2,3}^* + L_{2,4}^* \right) \right|
$$

$$
\leq \left| \mathbf{E}_{W,\epsilon} \left( \Psi \left( \frac{\sqrt{\chi_\alpha^2(1)} - L_{2,1}(L_{2,2} + L_{2,3}) + L_{2,4}}{1 + a_0(n)} \right) - \Psi \left( \sqrt{\chi_\alpha^2(1)} - L_{2,1}^* L_{2,3}^* + L_{2,4}^* \right) \right) \cdot \mathbf{1}_{\mathcal{A}} \right| \quad (58)
$$

$$
+ \left| \mathbf{E}_{W,\epsilon} \left( \Psi \left( \frac{\sqrt{\chi_\alpha^2(1)} - L_{2,1}(L_{2,2} + L_{2,3}) + L_{2,4}}{1 + a_0(n)} \right) - \Psi \left( \sqrt{\chi_\alpha^2(1)} - L_{2,1}^* L_{2,3}^* + L_{2,4}^* \right) \right) \cdot \mathbf{1}_{\mathcal{A}^c} \right|
$$

where $\mathcal{A}$ is defined in Lemma 1.

The following Lemma controls the terms in (58), whose proof is present in Section 5.2.

**Lemma 2.** *Under the same assumptions as Theorem 1,*

$$
\left| \mathbf{E}_{W,\epsilon} \left( \Psi \left( \frac{\sqrt{\chi_\alpha^2(1)} - L_{2,1}(L_{2,2} + L_{2,3}) + L_{2,4}}{1 + a_0(n)} \right) - \Psi \left( \sqrt{\chi_\alpha^2(1)} - L_{2,1}^* L_{2,3}^* + L_{2,4}^* \right) \right) \cdot \mathbf{1}_{\mathcal{A}^c} \right| \leq \mathbf{P}(\mathcal{A}^c), \quad (59)
$$

$$
\left| \mathbf{E}_{W,\epsilon} \left( \Psi \left( \frac{\sqrt{\chi_\alpha^2(1)} - L_{2,1}(L_{2,2} + L_{2,3}) + L_{2,4}}{1 + a_0(n)} \right) - \Psi \left( \sqrt{\chi_\alpha^2(1)} - L_{2,1}^* L_{2,3}^* + L_{2,4}^* \right) \right) \cdot \mathbf{1}_{\mathcal{A}} \right| \to 0. \quad (60)
$$

By (58), (59),(60) and Lemma 1, we establish (11) in the main paper.

### 4.3 Proof of Theorem 4

We start with introducing some notations. Define $\Pi = (\xi, \epsilon) \in \mathbb{R}^{n \times 2}$, $v = \sum_{j \in \mathcal{V}} \gamma_j \widehat{V}_{\cdot j}$, and

$$
\Delta^{\Theta_{11}} = \sqrt{n} \left( \widehat{\Theta}_{11} - \frac{1}{n} \Pi_{\cdot 1}' \Pi_{\cdot 1} \right), \quad \Delta^{\Theta_{12}} = \sqrt{n} \left( \widehat{\Theta}_{12} - \frac{1}{n} \Pi_{\cdot 1}' \Pi_{\cdot 2} \right),
$$

$$
\Delta^{\Theta_{22}} = \sqrt{n} \left( \widehat{\Theta}_{22} - \frac{1}{n} \Pi_{\cdot 2}' \Pi_{\cdot 2} \right).
$$

where $\widehat{V}_{\cdot j}, \widehat{\Theta}_{11}, \widehat{\Theta}_{22}$ and $\widehat{\Theta}_{12}$ are stated in Definition 2 of the well-behaved estimators. Let $P_v$ denote the projection matrix to the direction of $v$ and $P_{v^\perp}$ denote the projection matrix to the orthogonal complement of $v$, that is, $P_v = v(v'v)^{-1}v'$ and $P_{v^\perp} = \mathbf{I} - v(v'v)^{-1}v'$. To facilitate the discussion, we define the following events,

$$
\begin{aligned}
&\mathcal{B}_1 = \left\{\left|\frac{\Delta^\beta}{\sqrt{\text{Var}_1}}\right| \leq C\left(\sqrt{s_{z1}}\frac{s\log p}{\sqrt{n}} + \frac{1}{\|\gamma\|_2}\frac{s_{z1}\log p}{\sqrt{n}}\right)\right\}, \mathcal{B}_2 = \left\{\left|\widehat{\beta} - \beta\right| \leq C\sqrt{\frac{\log p}{n}}\sqrt{V_1}\right\}, \\
&\mathcal{B}_3 = \left\{\max|\widetilde{\gamma}_j - \gamma_j| \leq C\sqrt{\frac{\log p}{n}}\right\}, \mathcal{B}_4 = \left\{\max\{|\Delta^{\Theta_{12}}|,|\Delta^{\Theta_{22}}|,|\Delta^{\Theta_{11}}|\} \leq C\frac{s\log p}{\sqrt{n}}\right\}, \\
&\mathcal{B}_5 = \left\{\max_{1\leq i,j\leq 2}\left|\widehat{\Theta}_{ij} - \Theta_{ij}\right| \leq \sqrt{\frac{\log p}{n}}\right\}, \mathcal{B}_6 = \left\{c_0 \leq \min_{1\leq j\leq p_z}\|\widehat{V}_{\cdot j}\|_2 \leq \max_{1\leq j\leq p_z}\|\widehat{V}_{\cdot j}\|_2 \leq C_0\right\}, \\
&\mathcal{B}_7 = \left\{\|\sum_{j\in\mathcal{V}}\gamma_j\widehat{V}_{\cdot j}\|_2 \geq c_1\|\gamma\|_2\right\}, \mathcal{B}_8 = \left\{\frac{c_1^2}{\|\gamma_\mathcal{V}\|_2^2} \leq \frac{1}{\|\gamma_\mathcal{V}\|_2^4}\left\|\sum_{j\in\mathcal{V}}\gamma_j\widehat{V}_{\cdot j}\right\|_2^2 \leq \frac{C_0^2 s_{z1}}{\|\gamma_\mathcal{V}\|_2^2}\right\}, \\
&\mathcal{B}_9 = \left\{\frac{1}{\sqrt{n}}\left|\left(\Pi'_{\cdot 1}P_v\Pi_{\cdot 2} - \Theta_{12}\right) - \beta\left(\Pi'_{\cdot 2}P_v\Pi_{\cdot 2} - \Theta_{22}\right)\right| \leq C\sqrt{\frac{\log p}{n}}\right\}.
\end{aligned}
$$
(61)

and $\mathcal{B} = \cap_{i=1}^9 \mathcal{B}_i$. The following Lemma controls the probability of these events, whose proof is present in Section 5.4.

**Lemma 3.** *Under the same assumptions as Theorem 4, we have*

$$\liminf_{n\to\infty} \mathbf{P}(\mathcal{B}) = 1; \qquad (62)$$

The proof of the Theorem 4 depends on the following error decomposition of the proposed estimator $\widehat{\Sigma}_{12}$ stated in the following lemma, whose proof is present in Section 5.3.

**Lemma 4.** *Under the same assumptions as Theorem 4, then the following error decomposition holds,*

$$\sqrt{n}\frac{\widehat{\Sigma}_{12} - \Sigma_{12}}{\sqrt{\text{Var}_2 + \Theta_{22}^2\text{Var}_1}} = \frac{M_1 + M_2}{\sqrt{\text{Var}_2 + \Theta_{22}^2\text{Var}_1}} + \frac{R}{\sqrt{\text{Var}_2 + \Theta_{22}^2\text{Var}_1}}. \qquad (63)$$

$$M_1 \perp M_2 \mid W; \qquad (64)$$

where $M_1 = -\Theta_{22}\frac{1}{\sum_{j\in\mathcal{S}}\gamma_j^2}\sum_{j\in\mathcal{S}}\gamma_j(\widehat{V}_{\cdot j})'(\Pi_{\cdot 1} - \beta\Pi_{\cdot 2})$, $\text{Var}_1 = \Sigma_{11}\left\|\sum_{j\in\mathcal{S}}\gamma_j\widehat{V}_{\cdot j}/\sqrt{n}\right\|_2^2/\left(\sum_{j\in\mathcal{S}}\gamma_j^2\right)^2$, $M_2 = \frac{1}{\sqrt{n}}\left((\Pi'_{\cdot 1}P_{v^\perp}\Pi_{\cdot 2} - (n-1)\Theta_{12}) - \beta(\Pi'_{\cdot 2}P_{v^\perp}\Pi_{\cdot 2} - (n-1)\Theta_{22})\right)$, $\text{Var}_2 = \Theta_{11}\Theta_{22} +$

$\Theta_{12}^2 + 2\beta^2\Theta_{22}^2 - 4\beta\Theta_{12}\Theta_{22}$. In addition, on the event $\mathcal{B}$,

$$\left|\frac{R}{\sqrt{\text{Var}_2 + \Theta_{22}^2\text{Var}_1}}\right| \leq C\left(\sqrt{s_{z1}}\frac{s\log p}{\sqrt{n}} + \frac{1}{\|\gamma\|_2}\frac{s_{z1}\log p}{\sqrt{n}}\right), \tag{65}$$

and

$$\left|\sqrt{\frac{\widehat{\text{Var}}(\widehat{\Sigma}_{12})}{\Theta_{22}^2\text{Var}_1 + \text{Var}_2}} - 1\right| \leq C\frac{1}{\sqrt{s_{z1}\log p}}. \tag{66}$$

By the error decomposition of $\widehat{\Sigma}_{12}$ stated in Lemma 4, we establish the following expression for the power of $Q_E$,

$$\mathbf{P}\left(|Q_E| \geq z_{\alpha/2}\right) = \mathbf{P}\left(\left|\sqrt{n}\frac{\widehat{\Sigma}_{12} - \Sigma_{12}}{\sqrt{\text{Var}_2 + \Theta_{22}^2\text{Var}_1}} + \frac{\sqrt{n}\Sigma_{12}}{\sqrt{\Theta_{22}^2\text{Var}_1 + \text{Var}_2}}\right| \geq z_{\alpha/2}\sqrt{\frac{\widehat{\text{Var}}(\widehat{\Sigma}_{12})}{\Theta_{22}^2\text{Var}_1 + \text{Var}_2}}\right)$$

$$= \mathbf{P}\left(\frac{M_1 + M_2}{\sqrt{\text{Var}_2 + \Theta_{22}^2\text{Var}_1}} + B(\theta, W) \geq z_{\alpha/2}\sqrt{\frac{\widehat{\text{Var}}(\widehat{\Sigma}_{12})}{\Theta_{22}^2\text{Var}_1 + \text{Var}_2}} - \frac{R}{\sqrt{\text{Var}_2 + \Theta_{22}^2\text{Var}_1}}\right)$$

$$+ \mathbf{P}\left(\frac{M_1 + M_2}{\sqrt{\text{Var}_2 + \Theta_{22}^2\text{Var}_1}} + B(\theta, W) \leq -z_{\alpha/2}\sqrt{\frac{\widehat{\text{Var}}(\widehat{\Sigma}_{12})}{\Theta_{22}^2\text{Var}_1 + \text{Var}_2}} - \frac{R}{\sqrt{\text{Var}_2 + \Theta_{22}^2\text{Var}_1}}\right). \tag{67}$$

where $B(\theta, W) = \frac{\sqrt{n}\Sigma_{12}}{\sqrt{\Theta_{22}^2\text{Var}_1 + \text{Var}_2}}$.

The remaining proof is based on Lemma 5, which establishes that the power curve (67) converges to the power curve defined in (24). The essential idea of Lemma 5 is to establish the limiting distribution of the term $(M_1 + M_2)/\sqrt{\text{Var}_2 + \Theta_{22}^2\text{Var}_1}$ in (63) and show that the remainder $R/\sqrt{\text{Var}_2 + \Theta_{22}^2\text{Var}_1}$ is negligible compared to $(M_1+M_2)/\sqrt{\text{Var}_2 + \Theta_{22}^2\text{Var}_1}$. The proof of Lemma 5 can be found in Section 5.6.

**Lemma 5.** *Under the same assumptions as Theorem 4, we have*

$$\mathbf{P}\left(\frac{M_1 + M_2}{\sqrt{\text{Var}_2 + \Theta_{22}^2\text{Var}_1}} + B(\theta, W) \geq z_{\alpha/2}\sqrt{\frac{\widehat{\text{Var}}(\widehat{\Sigma}_{12})}{\Theta_{22}^2\text{Var}_1 + \text{Var}_2}} - \frac{R_1 + R_2 + R_3}{\sqrt{\text{Var}_2 + \Theta_{22}^2\text{Var}_1}}\right) \tag{68}$$
$$- \left(1 - \mathbf{E}_W\Phi\left(z_{\alpha/2} - B(\theta, W)\right)\right) \to 0.$$

*and*

$$\mathbf{P}\left(\frac{M_1 + M_2}{\sqrt{\text{Var}_2 + \Theta_{22}^2\text{Var}_1}} + B(\theta, W) \leq -z_{\alpha/2}\sqrt{\frac{\widehat{\text{Var}}(\widehat{\Sigma}_{12})}{\Theta_{22}^2\text{Var}_1 + \text{Var}_2}} - \frac{R_1 + R_2 + R_3}{\sqrt{\text{Var}_2 + \Theta_{22}^2\text{Var}_1}}\right) \tag{69}$$
$$- \mathbf{E}_W\Phi\left(-z_{\alpha/2} - B(\theta, W)\right) \to 0.$$

By (67) and Lemma 5, we establish (35). And (36) follows from (35) by taking $\Sigma_{12} = 0$.

## 4.4 Proof of Theorem 2

**Lemma 6.** *Under the same assumptions as Theorem 2, the same conclusions for Lemma 4 hold with*

$$M_1 = -\Theta_{22} \frac{1}{\sum_{j \in \mathcal{S}} \gamma_j^2} \sum_{j \in \mathcal{S}} \gamma_j (\widehat{V}_{\cdot j})' \left( \Pi_{\cdot 1} - \beta \Pi_{\cdot 2} \right),$$

$$M_2 = \frac{1}{\sqrt{n}} \left( \left( \Pi'_{\cdot 1} P_{v^\perp} \Pi_{\cdot 2} - (n-1)\Theta_{12} \right) - \beta \left( \Pi'_{\cdot 2} P_{v^\perp} \Pi_{\cdot 2} - (n-1)\Theta_{22} \right) \right),$$

$$\mathrm{Var}_1 = \Sigma_{11} \left\| \sum_{j \in \mathcal{S}} \gamma_j \widehat{V}_{\cdot j}/\sqrt{n} \right\|_2^2 \Big/ \left( \sum_{j \in \mathcal{S}} \gamma_j^2 \right)^2 \quad and \quad \mathrm{Var}_2 = \Theta_{11}\Theta_{22} + \Theta_{12}^2 + 2\beta^2 \Theta_{22}^2 - 4\beta\Theta_{12}\Theta_{22}.$$

Applying the similar argument with the proof of Theorem 4 in Section 4.3, we can establish Theorem 2.

# 5 Proof of key lemmas

## 5.1 Proof of Lemma 1

We first introduce the following technical lemmas, which will be used to prove Lemma 1. The first lemma (Theorem 2.3 in Boucheron et al. (2013)) is a concentration result of $\chi^2$ random variable.

**Lemma 7.** *Let $\chi_n^2$ denote the $\chi^2$ random variable with $n$ degrees of freedom, then we have the following concentration inequality,*

$$\mathbb{P}\left( \left| \chi_n^2 - E\chi_n^2 \right| > 2\sqrt{nt} + 2t \right) \leq 2\exp(-t).$$

The following lemma establishes the concentration of sum of independent centered subexponential random variables (Vershynin (2012) and Javanmard and Montanari (2014)),

**Lemma 8.** *Let $X_i$ denote Sub-exponential random variable with the Sub-exponential norm $K = \|X_i\|_{\psi_1}$, then we have*

$$P\left(\frac{1}{n}|\sum_{i=1}^{n} X_i| \geq \epsilon\right) \leq 2\exp\left(-\frac{1}{6}n\min\left(\frac{\epsilon}{K}, \frac{\epsilon^2}{K^2}\right)\right). \tag{70}$$

In the following, we first decompose some key terms $\frac{1}{n}D'P_{X^\perp}D$, $\frac{1}{n}D'(P_W - P_X)D$ and $\sqrt{(\frac{1}{n}D'(P_W - P_X)D)^{-1} - (\frac{1}{n}D'P_{X^\perp}D)^{-1}}$ in the definitions of $L_{2,i}$ for $i = 1, 3, 4$. Since $Z_{i\cdot} = X_{i\cdot}\left(\Lambda_{xx}^{-1}\Lambda_{xz}\right) + \bar{Z}_{i\cdot}$ where $\begin{pmatrix} \bar{Z}_{i\cdot} \\ \bar{X}_{i\cdot} \end{pmatrix} \sim N\left(\begin{pmatrix} 0 \\ 0 \end{pmatrix}, \begin{pmatrix} \Lambda_{\mathcal{I}|\mathcal{I}^c} & 0 \\ 0 & \Lambda_{xx} \end{pmatrix}\right)$, we have $D = \bar{Z}\gamma + X\left(\psi + \Lambda_{xx}^{-1}\Lambda_{xz}\gamma\right) + \epsilon$ and then further establish the following decompositions,

$$\frac{1}{n}D'P_{X^\perp}D = \frac{1}{n}\left(\bar{Z}\gamma + X\left(\psi + \Lambda_{xx}^{-1}\Lambda_{xz}\gamma\right) + \epsilon\right)'\left(I - X(X'X)^{-1}X'\right)\left(\bar{Z}\gamma + X\left(\psi + \Lambda_{xx}^{-1}\Lambda_{xz}\gamma\right) + \epsilon\right)$$
$$= \frac{1}{n}\left(\bar{Z}\gamma + \epsilon\right)'\left(I - X(X'X)^{-1}X'\right)\left(\bar{Z}\gamma + \epsilon\right),$$

$$\frac{1}{n}D'(P_W - P_X)D$$
$$= \frac{1}{n}\left(\bar{Z}\gamma + X\left(\psi + \Lambda_{xx}^{-1}\Lambda_{xz}\gamma\right) + \epsilon\right)'\left(W(W'W)^{-1}W' - X(X'X)^{-1}X'\right)\left(\bar{Z}\gamma + X\left(\psi + \Lambda_{xx}^{-1}\Lambda_{xz}\gamma\right) + \epsilon\right)$$
$$= \frac{1}{n}\left(\bar{Z}\gamma + \epsilon\right)'\left(W(W'W)^{-1}W' - X(X'X)^{-1}X'\right)\left(\bar{Z}\gamma + \epsilon\right)$$
$$= \frac{1}{n}\gamma'\bar{Z}'\left(I - X(X'X)^{-1}X'\right)\bar{Z}\gamma + \frac{1}{n}\epsilon'\left(W(W'W)^{-1}W' - X(X'X)^{-1}X'\right)\epsilon + \frac{1}{n}\epsilon'\left(I - X(X'X)^{-1}X'\right)\bar{Z}\gamma,$$

and

$$\sqrt{(\frac{1}{n}D'(P_W - P_X)D)^{-1} - (\frac{1}{n}D'P_{X^\perp}D)^{-1}} = \sqrt{\frac{\frac{1}{n}D'\left(P_{X^\perp} - (P_W - P_X)\right)D}{(\frac{1}{n}D'(P_W - P_X)D)(\frac{1}{n}D'P_{X^\perp}D)}}. \tag{71}$$

The proof also relies on the following expressions $\frac{1}{n}D'\left(P_{X^\perp} - (P_W - P_X)\right)\epsilon = \frac{1}{n}\epsilon'\left(I - W(W'W)^{-1}W'\right)\epsilon$, and $\frac{1}{n}D'\left(P_{X^\perp} - (P_W - P_X)\right)D = \frac{1}{n}\epsilon'\left(I - W(W'W)^{-1}W'\right)\epsilon$. We introduce the following quantities to represent the approximation errors,

$$\omega_1 = \omega_1(n) = \frac{\frac{1}{n}\left(\bar{Z}\gamma + \epsilon\right)'\left(I - X\left(X'X\right)^{-1}X'\right)\left(\bar{Z}\gamma + \epsilon\right)}{\frac{n-p_x}{n}\left(\gamma'\Lambda_{\mathcal{I}|\mathcal{I}^c}\gamma + \Sigma_{22}\right)} - 1$$

$$\omega_2 = \omega_2(n) = \frac{\frac{1}{n}\gamma'\bar{Z}'\left(I - X\left(X'X\right)^{-1}X'\right)\bar{Z}\gamma}{\frac{n-p_x}{n}\gamma'\Lambda_{\mathcal{I}|\mathcal{I}^c}\gamma} - 1$$

$$\omega_3 = \omega_3(n) = \frac{\frac{1}{n}\epsilon'\left(W\left(W'W\right)^{-1}W' - X\left(X'X\right)^{-1}X'\right)\epsilon}{\frac{p_z}{n}\Sigma_{22}} - 1 \tag{72}$$

$$\omega_4 = \omega_4(n) = \frac{\frac{1}{n}\epsilon'\left(I - W\left(W'W\right)^{-1}W'\right)\epsilon}{\frac{n-p}{n}\Sigma_{22}} - 1$$

$$\omega_5 = \omega_5(n) = \frac{1}{n}\gamma'\bar{Z}'\left(I - X\left(X'X\right)^{-1}X'\right)\epsilon.$$

Define the following events

$$\mathcal{A}_1 = \left\{|\omega_1| \le 2\sqrt{\frac{\log(n - p_x)}{n - p_x}} + 2\frac{\log(n - p_x)}{n - p_x}\right\}$$

$$\mathcal{A}_2 = \left\{|\omega_2| \le 2\sqrt{\frac{\log(n - p_x)}{n - p_x}} + 2\frac{\log(n - p_x)}{n - p_x}\right\}$$

$$\mathcal{A}_3 = \left\{|\omega_3| \le 2\sqrt{\frac{\log(p_z)}{p_z}} + 2\frac{\log(p_z)}{p_z}\right\} \tag{73}$$

$$\mathcal{A}_4 = \left\{|\omega_4| \le 2\sqrt{\frac{\log(n - p)}{n - p}} + 2\frac{\log(n - p)}{n - p}\right\}$$

$$\mathcal{A}_5 = \left\{|\omega_5| \le C\frac{\sqrt{(n - p_x)\log(n - p_x)}}{n}\sqrt{\Sigma_{22} \cdot \gamma'\Lambda_{\mathcal{I}|\mathcal{I}^c}\gamma}\right\}$$

and define

$$\mathcal{A} = \cap_{i=1}^{5}\mathcal{A}_i. \tag{74}$$

In the following, we will first control the probability of event $\mathcal{A}$ defined in (74) and then show that (56) and (57) holds on the event $\mathcal{A}$.

<u>Proof of (55)</u> Conditioning on $X$, we have

$$\omega_1 + 1 \sim \frac{1}{n - p_x}\chi^2_{n-p_x} \quad \text{and} \quad \omega_2 + 1 \sim \frac{1}{n - p_x}\chi^2_{n-p_x}.$$

Conditioning on $W$, we have

$$\omega_3 + 1 \sim \frac{1}{p_z}\chi^2_{p_z} \quad \text{and} \quad \omega_4 + 1 \sim \frac{1}{n-p}\chi^2_{n-p}.$$

By Lemma 7, we can establish that

$$\mathbf{P}\left(\mathcal{A}_1 \cap \mathcal{A}_2 \cap \mathcal{A}_3 \cap \mathcal{A}_4\right) \geq 1 - (\min\{n-p, p_z\})^{-c}. \tag{75}$$

Conditioning on $X$, $\omega_5 = \frac{1}{n}\sum_{i=1}^{n-p_x} U_i V_i$, where $U_i$ is independent of $V_i$, $U_i$ follows i.i.d normal with mean 0 and variance $\gamma'\Lambda_{\mathcal{I}|\mathcal{I}^c}\gamma$ and $V_i$ follows i.i.d normal with mean 0 and variance $\Sigma_{22}$. Note that $K = \|U_i V_i\|_{\psi_1} \leq 2\sqrt{\gamma'\Lambda_{\mathcal{I}|\mathcal{I}^c}\gamma}\sqrt{\Sigma_{22}}$. By Lemma 8, we establish that

$$\mathbf{P}\left(\mathcal{A}_5\right) \geq 1 - (n-p_x)^{-c}. \tag{76}$$

The control of probability (55) follows from (75) and (76).

Proof of (56) and (57) By the definition (72), we have the following expressions,

$$\begin{aligned}
\frac{1}{n}\epsilon'\left(W\left(W'W\right)^{-1}W' - X\left(X'X\right)^{-1}X'\right)\epsilon &= \frac{p_z}{n}\Sigma_{22}(1+\omega_3), \\
\frac{1}{n}D'P_{X^\perp}D &= \frac{n-p_x}{n}\left(\gamma'\Lambda_{\mathcal{I}|\mathcal{I}^c}\gamma + \Sigma_{22}\right)(1+\omega_1), \\
\frac{1}{n}D'(P_W - P_X)D &= \frac{n-p_x}{n}\left(\gamma'\Lambda_{\mathcal{I}|\mathcal{I}^c}\gamma\right)(1+\omega_2) + \frac{p_z}{n}\Sigma_{22}(1+\omega_3) + \omega_5, \\
\frac{1}{n}D'\left(P_{X^\perp} - (P_W - P_X)\right)\epsilon &= \frac{n-p}{n}\Sigma_{22}(1+\omega_4), \\
\frac{1}{n}D'\left(P_{X^\perp} - (P_W - P_X)\right)D &= \frac{n-p}{n}\Sigma_{22}(1+\omega_4).
\end{aligned} \tag{77}$$

By the first equation of (77), we have

$$L_{2,3} = L^*_{2,3}(1+\omega_4). \tag{78}$$

Define

$$h_1 = h_1 = \frac{\frac{1}{n}D'(P_W - P_X)D}{\left(\frac{n-p_x}{n}\left(\gamma'\Lambda_{\mathcal{I}|\mathcal{I}^c}\gamma\right) + \frac{p_z}{n}\Sigma_{22}\right)} - 1. \tag{79}$$

Note that

$$\frac{1}{n}D'(P_W - P_X)D = \frac{n-p_x}{n}\left(\gamma'\Lambda_{\mathcal{I}|\mathcal{I}^c}\gamma\right)\left(1 + \omega_2 + \frac{\omega_5}{\frac{n-p_x}{n}\left(\gamma'\Lambda_{\mathcal{I}|\mathcal{I}^c}\gamma\right)}\right) + \frac{p_z}{n}\Sigma_{22}\left(1+\omega_3\right).$$

Hence, on the event $\mathcal{A}$,

$$|h_1| \leq |\omega_2| + \left|\frac{\omega_5}{\frac{n-p_x}{n}\left(\gamma'\Lambda_{\mathcal{I}|\mathcal{I}^c}\gamma\right)}\right| + |\omega_3| \leq C\left(\frac{\log p_z}{p_z} + \sqrt{\frac{\log(n-p_x)}{n-p_x}}\sqrt{\frac{\Sigma_{22}}{\gamma'\Lambda_{\mathcal{I}|\mathcal{I}^c}\gamma}} + \frac{\log(n-p_x)}{n-p_x}\right).$$

By plugging the second and fifth equation of (72) and (79) into (71), we have the following key expressions,

$$\sqrt{(\frac{1}{n}D'(P_W - P_X)D)^{-1} - (\frac{1}{n}D'P_{X^\perp}D)^{-1}}$$
$$= \sqrt{\frac{\frac{n-p}{n}\Sigma_{22}}{\left(\frac{n-p_x}{n}\left(\gamma'\Lambda_{\mathcal{I}|\mathcal{I}^c}\gamma + \Sigma_{22}\right)\right)\left(\frac{n-p_x}{n}\left(\gamma'\Lambda_{\mathcal{I}|\mathcal{I}^c}\gamma\right) + \frac{p_z}{n}\Sigma_{22}\right)}}\sqrt{\frac{1+\omega_4}{(1+\omega_1)(1+h_1)}}$$

and hence

$$L_{2,1} = L_{2,1}^*\sqrt{\frac{1+\omega_4}{(1+\omega_1)(1+h_1)}} \tag{80}$$

Combined with the second and forth equation of (77), we have

$$\frac{(\frac{1}{n}D'P_{X^\perp}D)^{-1}\frac{1}{n}D'\left(P_{X^\perp} - (P_W - P_X)\right)\epsilon}{\sqrt{(\frac{1}{n}D'(P_W - P_X)D)^{-1} - (\frac{1}{n}D'P_{X^\perp}D)^{-1}}} = \sqrt{\frac{(1-\frac{p}{n})\Sigma_{22}\left(\frac{\gamma'\Lambda_{\mathcal{I}|\mathcal{I}^c}\gamma}{p_z\Sigma_{22}} + \frac{1}{n-p_x}\right)}{\frac{\gamma'\Lambda_{\mathcal{I}|\mathcal{I}^c}\gamma}{p_z\Sigma_{22}} + \frac{1}{p_z}}}\sqrt{\frac{(1+\omega_1)(1+h_1)}{1+\omega_4}}$$

and hence

$$L_{2,4} = L_{2,4}^*\sqrt{\frac{(1+\omega_1)(1+h_1)}{1+\omega_4}} \tag{81}$$

By defining $\mathcal{A}$ as in (74), the control of terms (56) and (57) follows from (80), (78), (81) and (73).

## 5.2 Proof of Lemma 2

The proof of (59) follows from $\sup_{x \in \mathbb{R}} |\Psi(x)| \leq 1$. It remains to prove (60). By the fact that $\sup_{x \in \mathbb{R}} |\Psi'(x)| \leq 1$, we have

$$\left| \mathbf{E}_{W,\epsilon} \left( \Psi \left( \frac{\sqrt{\chi_\alpha^2(1)} - L_{2,1}(L_{2,2} + L_{2,3}) + L_{2,4}}{1 + a_0(n)} \right) - \Psi \left( \sqrt{\chi_\alpha^2(1)} - L_{2,1}^* L_{2,3}^* + L_{2,4}^* \right) \right) \cdot \mathbf{1}_{\mathcal{A}} \right|$$

$$\leq \mathbf{E}_{W,\epsilon} \left| \frac{\sqrt{\chi_\alpha^2(1)} - L_{2,1}(L_{2,2} + L_{2,3}) + L_{2,4}}{1 + a_0(n)} - \left( \sqrt{\chi_\alpha^2(1)} - L_{2,1}^* L_{2,3}^* + L_{2,4}^* \right) \right| \cdot \mathbf{1}_{\mathcal{A}}$$

$$\leq \left| \frac{a_0(n)}{1 + a_0(n)} \right| \sqrt{\chi_\alpha^2(1)} + |L_{2,1}^*| \times \mathbf{E}_{W,\epsilon} \left( \left| 1 - \frac{L_{2,1}}{L_{2,1}^*(1 + a_0(n))} \right| + 1 \right) |L_{2,2}| \cdot \mathbf{1}_{\mathcal{A}}$$

$$+ |L_{2,1}^* L_{2,3}^*| \times \mathbf{E}_{W,\epsilon} \left| 1 - \frac{L_{2,1} L_{2,3}}{L_{2,1}^* L_{2,3}^*(1 + a_0(n))} \right| \cdot \mathbf{1}_{\mathcal{A}} + |L_{2,4}^*| \times \mathbf{E}_{W,\epsilon} \left| 1 - \frac{L_{2,4}}{L_{2,4}^*(1 + a_0(n))} \right| \cdot \mathbf{1}_{\mathcal{A}}$$

(82)

where the last inequality follows from triangle inequality.

By the definition (54), the last term in the above inequality can be expressed as

$$\left| \frac{a_0(n)}{1 + a_0(n)} \right| \sqrt{\chi_\alpha^2(1)} + |L_{2,1}^*| \times \mathbf{E}_{W,\epsilon} \left( \left| 1 - \frac{1 + a_1(n)}{1 + a_0(n)} \right| + 1 \right) \cdot |L_{2,2}| \cdot \mathbf{1}_{\mathcal{A}}$$

$$+ |L_{2,1}^* L_{2,3}^*| \times \mathbf{E}_{W,\epsilon} \left| 1 - \frac{(1 + a_1(n))(1 + a_2(n))}{1 + a_0(n)} \right| \cdot \mathbf{1}_{\mathcal{A}} + |L_{2,4}^*| \times \mathbf{E}_{W,\epsilon} \left| 1 - \frac{1 + a_3(n)}{1 + a_0(n)} \right| \cdot \mathbf{1}_{\mathcal{A}},$$

$$\leq \frac{C}{n} + C \left( \frac{\log p_z}{p_z} + \frac{\log(n-p)}{n-p} + \sqrt{\frac{\log(n-p_x)}{n-p_x}} \sqrt{\frac{\Sigma_{22}}{\gamma' \Lambda_{\mathcal{I} | \mathcal{I}^c} \gamma}} \right)$$

$$\left( |L_{2,1}^*| \mathbf{E}_{W,\epsilon} |L_{2,2}| + |L_{2,1}^*| \cdot |L_{2,3}^*| + |L_{2,4}^*| \right),$$

(83)

which the last inequality follows the fact $|a_0(n)| \leq C \frac{1}{n}$ and Lemma 1. By the parametrization $\Sigma_{12} = \frac{\Delta_1}{\sqrt{n}}$, we have $|L_{2,4}^*| \leq \frac{\Delta_1}{\sqrt{\Sigma_{11} \Sigma_{22}}}$ and

$$|L_{2,1}^*| \cdot |L_{2,3}^*| = \frac{\Delta_1}{\sqrt{\Sigma_{11} \Sigma_{22}}} \sqrt{\frac{\frac{n-p}{(n-p_x)^2 n}}{\left( \frac{\gamma' \Lambda_{\mathcal{I} | \mathcal{I}^c} \gamma}{p_z \Sigma_{22}} + \frac{1}{n-p_x} \right) \left( \frac{\gamma' \Lambda_{\mathcal{I} | \mathcal{I}^c} \gamma}{p_z \Sigma_{22}} + \frac{1}{p_z} \right)}}.$$

It remains to control $|L_{2,1}^*| \mathbf{E}_{W,\epsilon} |L_{2,2}|$. For the term $L_{2,2}$, conditioning on $X$,

$$L_{2,2} = \frac{1}{n} \sum_{i=1}^{n-p_x} U_i V_i = \frac{p_z}{n} \Sigma_{22} \sqrt{\frac{n-p_x}{p_z}} \sqrt{\frac{\gamma' \Lambda_{\mathcal{I} | \mathcal{I}^c} \gamma}{p_z \Sigma_{22}}} \left( \frac{1}{\sqrt{n-p_x}} \sum_{i=1}^{n-p_x} \frac{U_i}{\sqrt{\gamma' \Lambda_{\mathcal{I} | \mathcal{I}^c} \gamma}} \frac{V_i}{\sqrt{\Sigma_{22}}} \right), \quad (84)$$

where $U_i$ is independent of $V_i$, $U_i$ follows i.i.d normal with mean 0 and variance $\gamma' \Lambda_{\mathcal{I} | \mathcal{I}^c} \gamma$ and $V_i$

follows i.i.d normal with mean 0 and variance $\Sigma_{22}$. Based on the decomposition (84), we have

$$\left|L_{2,1}^*\right| \mathbf{E}_{W,\epsilon} |L_{2,2}| = \left|L_{2,1}^*\right| \cdot \left|L_{2,3}^*\right| \mathbf{E} \sqrt{\frac{n-p_{\mathrm{x}}}{p_{\mathrm{z}}}} \sqrt{\frac{\gamma'\Lambda_{\mathcal{I}|\mathcal{I}^c}\gamma}{p_{\mathrm{z}}\Sigma_{22}}} \left|\frac{1}{\sqrt{n-p_{\mathrm{x}}}} \sum_{i=1}^{n-p_{\mathrm{x}}} \frac{U_i}{\sqrt{\gamma'\Lambda_{\mathcal{I}|\mathcal{I}^c}\gamma}} \frac{V_i}{\sqrt{\Sigma_{22}}}\right|$$

$$\leq C \frac{\Delta_1}{\sqrt{\Sigma_{11}\Sigma_{22}}} \sqrt{\frac{\frac{n-p}{(n-p_{\mathrm{x}})p_{\mathrm{z}}n} \frac{\gamma'\Lambda_{\mathcal{I}|\mathcal{I}^c}\gamma}{p_{\mathrm{z}}\Sigma_{22}}}{\left(\frac{\gamma'\Lambda_{\mathcal{I}|\mathcal{I}^c}\gamma}{p_{\mathrm{z}}\Sigma_{22}} + \frac{1}{n-p_{\mathrm{x}}}\right)\left(\frac{\gamma'\Lambda_{\mathcal{I}|\mathcal{I}^c}\gamma}{p_{\mathrm{z}}\Sigma_{22}} + \frac{1}{p_{\mathrm{z}}}\right)}} \leq C \frac{\Delta_1}{\sqrt{\Sigma_{11}\Sigma_{22}}} \sqrt{\frac{\frac{n-p}{(n-p_{\mathrm{x}})p_{\mathrm{z}}n}}{\left(\frac{\gamma'\Lambda_{\mathcal{I}|\mathcal{I}^c}\gamma}{p_{\mathrm{z}}\Sigma_{22}} + \frac{1}{p_{\mathrm{z}}}\right)}}$$

By the fact that $(n-p_{\mathrm{x}})n\left(\frac{\gamma'\Lambda_{\mathcal{I}|\mathcal{I}^c}\gamma}{p_{\mathrm{z}}\Sigma_{22}} + \frac{1}{n-p_{\mathrm{x}}}\right)\left(\frac{\gamma'\Lambda_{\mathcal{I}|\mathcal{I}^c}\gamma}{p_{\mathrm{z}}\Sigma_{22}} + \frac{1}{p_{\mathrm{z}}}\right) \to \infty$ and $p_{\mathrm{z}}n\left(\frac{\gamma'\Lambda_{\mathcal{I}|\mathcal{I}^c}\gamma}{p_{\mathrm{z}}\Sigma_{22}} + \frac{1}{p_{\mathrm{z}}}\right) \to \infty$, we have

$$\left|L_{2,1}^*\right| \cdot \left|L_{2,3}^*\right| \to 0 \quad \text{and} \quad \left|L_{2,1}^*\right| \mathbf{E}_{W,\epsilon} |L_{2,2}| \to 0. \tag{85}$$

Since $\sqrt{C(\mathcal{I})} \gg \sqrt{\log(n-p_{\mathrm{x}})/(n-p_{\mathrm{x}})p_{\mathrm{z}}}$, combined with (82), (83) and (85), we establish (60).

### 5.3 Proof of Lemma 4

The error decomposition of $\widehat{\Sigma}_{12}$ depends on the error decomposition of $\widehat{\beta}$ defined in (33), which was obtained as Theorem 2 in Guo et al. (2016a) and stated in the following lemma.

**Lemma 9.** *Under the same assumptions as Theorem 4, the following property holds for the estimator $\widehat{\beta}$ defined in (33),*

$$\sqrt{n}\left(\widehat{\beta} - \beta\right) = T^\beta + \Delta^\beta \tag{86}$$

*with $T^\beta = \frac{1}{\sum_{j \in \mathcal{V}} \gamma_j^2} \sum_{j \in \mathcal{V}} \gamma_j (\widehat{V}_{.j})' (\Pi_{.1} - \beta \Pi_{.2})$ and*

$$\lim_{n \to \infty} \mathbf{P}\left(\left|\frac{\Delta^\beta}{\sqrt{\mathrm{Var}_1}}\right| \geq C\left(\sqrt{s_{z1}}\frac{s\log p}{\sqrt{n}} + \frac{1}{\|\gamma\|_2}\frac{s_{z1}\log p}{\sqrt{n}}\right)\right) = 0, \tag{87}$$

*where $\mathrm{Var}_1 = 1/\left(\sum_{j \in \mathcal{V}} \gamma_j^2\right)^2 \times \left\|\sum_{j \in \mathcal{V}} \gamma_j \widehat{V}_{.j}\right\|_2^2 (\Theta_{11} + \beta^2 \Theta_{22} - 2\beta\Theta_{12})$. Note that $T^\beta \mid W \sim N(0, \mathrm{Var}_1)$.*

By the definition $\widehat{\Sigma}_{12} = \widehat{\Theta}_{12} - \widehat{\beta}\widehat{\Theta}_{22}$, we have the following expression for $\widehat{\Sigma}_{12} - \Sigma_{12}$,

$$\widehat{\Sigma}_{12} - \Sigma_{12} = \left(\widehat{\Theta}_{12} - \widehat{\beta}\widehat{\Theta}_{22}\right) - (\Theta_{12} - \beta\Theta_{22})$$
$$= \left(\widehat{\Theta}_{12} - \Theta_{12}\right) - \beta\left(\widehat{\Theta}_{22} - \Theta_{22}\right) - \left(\widehat{\beta} - \beta\right)\Theta_{22} - \left(\widehat{\beta} - \beta\right)\left(\widehat{\Theta}_{22} - \Theta_{22}\right). \quad (88)$$

<u>Proof of (63)</u> By plugging the error bound for $\widehat{\beta} - \beta$ in Lemma 9 and the following error bounds for $(\widehat{\Theta}_{12} - \Theta_{12})$ and $(\widehat{\Theta}_{22} - \Theta_{22})$ into (88),

$$\sqrt{n}\left(\widehat{\Theta}_{12} - \Theta_{12}\right) = \frac{1}{\sqrt{n}}\left(\left(\Pi'_{\cdot 1}P_{v^\perp}\Pi_{\cdot 2} - (n-1)\Theta_{12}\right) + \left(\Pi'_{\cdot 1}P_v\Pi_{\cdot 2} - \Theta_{12}\right)\right) + \Delta^{\Theta_{12}},$$

$$\sqrt{n}\left(\widehat{\Theta}_{22} - \Theta_{22}\right) = \frac{1}{\sqrt{n}}\left(\left(\Pi'_{\cdot 2}P_{v^\perp}\Pi_{\cdot 2} - (n-1)\Theta_{22}\right) + \left(\Pi'_{\cdot 2}P_v\Pi_{\cdot 2} - \Theta_{22}\right)\right) + \Delta^{\Theta_{22}},$$

we establish (63) with the remainder term $R$ in (63) decomposed as $R = R_1 + R_2 + R_3$ where $R_1 = -\Theta_{22}\Delta^\beta + \Delta^{\Theta_{12}} - \beta\Delta^{\Theta_{22}}$, $R_2 = \frac{1}{\sqrt{n}}\left(\left(\Pi'_{\cdot 1}P_v\Pi_{\cdot 2} - \Theta_{12}\right) - \beta\left(\Pi'_{\cdot 2}P_v\Pi_{\cdot 2} - \Theta_{22}\right)\right)$ and $R_3 = \sqrt{n}\left(\widehat{\beta} - \beta\right)\left(\widehat{\Theta}_{22} - \Theta_{22}\right)$.

<u>Proof of (64)</u> Conditioning on $W$, $\begin{pmatrix} v' & 0 \\ 0 & v' \\ P_{v^\perp} & 0 \\ 0 & P_{v^\perp} \end{pmatrix}\begin{pmatrix} \Pi_{\cdot 1} \\ \Pi_{\cdot 2} \end{pmatrix}$ is jointly normal distribution. Since

$$\text{Cov}\left(\begin{pmatrix} v'\Pi_{\cdot 1} \\ v'\Pi_{\cdot 2} \end{pmatrix}, \begin{pmatrix} P_{v^\perp}\Pi_{\cdot 1} \\ P_{v^\perp}\Pi_{\cdot 2} \end{pmatrix} \mid W\right) = \mathbf{0},$$ we establish that conditioning on $W$,

$$\begin{pmatrix} \frac{1}{\|v\|_2}v'\Pi_{\cdot 1} \\ \frac{1}{\|v\|_2}v'\Pi_{\cdot 2} \end{pmatrix} \perp \begin{pmatrix} \frac{1}{\sqrt{n}}\left(\Pi'_{\cdot 1}P_{v^\perp}\Pi_{\cdot 2} - (n-1)\Theta_{12}\right) \\ \frac{1}{\sqrt{n}}\left(\Pi'_{\cdot 2}P_{v^\perp}\Pi_{\cdot 2} - (n-1)\Theta_{22}\right) \end{pmatrix},$$

and hence establish (64).

<u>Proof of (65)</u> On $\mathcal{B}_1 \cap \mathcal{B}_4$, we have

$$\left|\frac{R_1}{\sqrt{\text{Var}_2 + \Theta_{22}^2\text{Var}_1}}\right| \leq C\left(\sqrt{s_{z1}}\frac{s\log p}{\sqrt{n}} + \frac{1}{\|\gamma\|_2}\frac{s_{z1}\log p}{\sqrt{n}}\right). \quad (89)$$

On $\mathcal{B}_8$, we have
$$\left|\frac{R_2}{\sqrt{\text{Var}_2 + \Theta_{22}^2\text{Var}_1}}\right| \le C\sqrt{\frac{\log p}{n}}. \tag{90}$$

On $\mathcal{B}_2 \cap \mathcal{B}_5$, we have
$$\left|\frac{R_3}{\sqrt{\text{Var}_2 + \Theta_{22}^2\text{Var}_1}}\right| \le C\frac{\log p}{\sqrt{n}}. \tag{91}$$

By (89), (90) and (91), we establish (65).

<u>Proof of (66)</u> By the decomposition,

$$\sqrt{\frac{\widehat{\text{Var}}(\widehat{\Sigma}_{12})}{\Theta_{22}^2\text{Var}_1 + \text{Var}_2}} - 1 = \frac{\widehat{\Theta}_{22}^2\widehat{\text{Var}}_1 + \widehat{\text{Var}}_2 - \Theta_{22}^2\text{Var}_1 - \text{Var}_2}{\sqrt{\Theta_{22}^2\text{Var}_1 + \text{Var}_2}\left(\sqrt{\widehat{\Theta}_{22}^2\widehat{\text{Var}}_1 + \widehat{\text{Var}}_2} + \sqrt{\Theta_{22}^2\text{Var}_1 + \text{Var}_2}\right)},$$

we have

$$\left|\sqrt{\frac{\widehat{\text{Var}}(\widehat{\Sigma}_{12})}{\Theta_{22}^2\text{Var}_1 + \text{Var}_2}} - 1\right| \le \frac{\left|\widehat{\Theta}_{22}^2 - \Theta_{22}^2\right|\text{Var}_1 + \left|\widehat{\Theta}_{22}^2 - \Theta_{22}^2\right|\left|\widehat{\text{Var}}_1 - \text{Var}_1\right|}{\Theta_{22}^2\text{Var}_1 + \text{Var}_2}$$
$$+ \frac{\Theta_{22}^2\left|\widehat{\text{Var}}_1 - \text{Var}_1\right| + \left|\widehat{\text{Var}}_2 - \text{Var}_2\right|}{\Theta_{22}^2\text{Var}_1 + \text{Var}_2} \tag{92}$$

Note that
$$\left|\widehat{\text{Var}}_2 - \text{Var}_2\right| \le C\max\left\{\max_{1\le i,j\le 2}\left|\widehat{\Theta}_{ij} - \Theta_{ij}\right|, \left|\widehat{\beta} - \beta\right|\right\}$$

and hence
$$\frac{\left|\widehat{\text{Var}}_2 - \text{Var}_2\right|}{\Theta_{22}^2\text{Var}_1 + \text{Var}_2} \le C\sqrt{\frac{\log p}{n}}. \tag{93}$$

Also,
$$\frac{\left|\widehat{\Theta}_{22}^2 - \Theta_{22}^2\right|\text{Var}_1}{\Theta_{22}^2\text{Var}_1 + \text{Var}_2} \le C\sqrt{\frac{\log p}{n}}. \tag{94}$$

To establish (66), it remains to control $\frac{|\widehat{\text{Var}}_1 - \text{Var}_1|}{\Theta_{22}^2\text{Var}_1 + \text{Var}_2}$, which is further upper bounded by

$$\left|\frac{\widehat{\text{Var}}_1}{\text{Var}_1} - 1\right| = \left|\frac{\|\widetilde{\gamma}_{\mathcal{V}}\|_2^2}{\|\gamma\|_2^2}\frac{\left\|\sum_{j\in\widehat{\mathcal{V}}}\widetilde{\gamma}_j\widehat{V}_{\cdot j}\right\|_2^2}{\left\|\sum_{j\in\mathcal{V}}\gamma_j\widehat{V}_{\cdot j}\right\|_2^2}\frac{\widehat{\Theta}_{11} + \widehat{\beta}^2\widehat{\Theta}_{22} - 2\widehat{\beta}\widehat{\Theta}_{12}}{\Theta_{11} + \beta^2\Theta_{22} - 2\beta\Theta_{12}} - 1\right| \tag{95}$$

By (149) and (150) in Guo et al. (2016b), we have

$$\left|\frac{\widetilde{\|\gamma\|_2^2}}{\|\gamma_\mathcal{V}\|_2^2} - 1\right| \leq C \frac{1}{\|\gamma_\mathcal{V}\|_2^2} \left(s_{z1} \frac{\log p_z}{n} + C s_{z1} \left(s \frac{\log p}{n}\right)^2 + C\|\gamma_\mathcal{V}\|_2 \sqrt{\frac{2 s_{z1} \log p_z}{n}}\right), \quad (96)$$

and

$$\left|\frac{\left\|\sum_{j \in \widehat{\mathcal{V}}} \widetilde{\gamma}_j \widehat{V}_{\cdot j}\right\|_2}{\left\|\sum_{j \in \mathcal{V}} \gamma_j \widehat{V}_{\cdot j}\right\|_2} - 1\right| \leq C s_{z1} \sqrt{\frac{\log p}{n}}. \quad (97)$$

Note that

$$\left|\frac{\widehat{\Theta}_{11} + \widehat{\beta}^2 \widehat{\Theta}_{22} - 2\widehat{\beta}\widehat{\Theta}_{12}}{\Theta_{11} + \beta^2 \Theta_{22} - 2\beta \Theta_{12}} - 1\right| \leq C \max\left\{\max_{1 \leq i,j \leq 2}\left|\widehat{\Theta}_{ij} - \Theta_{ij}\right|, \left|\widehat{\beta} - \beta\right|\right\} \leq C \sqrt{\frac{\log p}{n}} \left(1 + \frac{\sqrt{s_{z1}}}{\|\gamma_\mathcal{V}\|_2}\right), \quad (98)$$

where the last inequality follows from the definition of $\mathcal{B}_2$ and $\mathcal{B}_5$ and the fact that $\sqrt{\text{Var}_1} \leq C \frac{\sqrt{s_{z1}}}{\|\gamma_\mathcal{V}\|_2}$. Combining (95), (96), (97) and (98), we establish that

$$\left|\frac{\widehat{\text{Var}}_1}{\text{Var}_1} - 1\right| \leq C \frac{1}{\|\gamma_\mathcal{V}\|_2^2} \left(s_{z1} \frac{\log p_z}{n} + C s_{z1} \left(s \frac{\log p}{n}\right)^2 + C\|\gamma_\mathcal{V}\|_2 \sqrt{\frac{2 s_{z1} \log p_z}{n}}\right)$$
$$+ C s_{z1} \sqrt{\frac{\log p}{n}} + C \sqrt{\frac{\log p}{n}} \left(1 + \frac{\sqrt{s_{z1}}}{\|\gamma_\mathcal{V}\|_2}\right) \leq C \frac{1}{\sqrt{s_{z1} \log p}}, \quad (99)$$

where the second inequality follows from the assumption $\|\gamma_\mathcal{V}^*\|_2 \gg s_{z1} \log p / \sqrt{n}$ and $\frac{\sqrt{s_{z1}} s \log p}{n} \to 0$. Combing (93), (94) and (99), we establish (66).

## 5.4 Proof of Lemma 3

The probability control of $\mathcal{B}$ in (62) in Lemma 3 follows from Lemma 9, the fact $\mathcal{B}_6 \cap \mathcal{B}_7 \subset \mathcal{B}_8$ and the fact that $\left(\widetilde{\gamma}, \widetilde{\Gamma}, \widehat{\Theta}_{11}, \widehat{\Theta}_{22}, \widehat{\Theta}_{12}\right)$ is well-behaved estimator,

## 5.5 Proof of Lemma 6

The proof of Lemma 6 is similar to that of Lemma 4 in Section 5.3. The major change is to replace Lemma 9 in Section 5.3 with the following lemma, which was Theorem 3 in Guo et al. (2016a).

**Lemma 10.** *Under the same assumptions as Theorem 2, the following property holds for the estimator $\widehat{\beta}$ defined in (22) in the main paper,*

$$\sqrt{n}\left(\widehat{\beta} - \beta\right) = T^{\beta} + \Delta^{\beta}, \tag{100}$$

*with $T^{\beta} = \frac{1}{\sum_{j \in \mathcal{S}} \gamma_j^2} \sum_{j \in \mathcal{S}} \gamma_j (\widehat{V}_{\cdot j})' (\Pi_{\cdot 1} - \beta \Pi_{\cdot 2})$ and*

$$\limsup_{n \to \infty} \mathbf{P}\left(\left|\frac{\Delta^{\beta}}{\sqrt{\operatorname{Var}_1}}\right| \geq C\left(\sqrt{s_{z1}} \frac{s \log p}{\sqrt{n}} + \frac{1}{\|\gamma\|_2} \frac{s_{z1} \log p}{\sqrt{n}}\right)\right) = 0, \tag{101}$$

*where $\operatorname{Var}_1 = 1/\left(\sum_{j \in \mathcal{S}} \gamma_j^2\right)^2 \times \left\|\sum_{j \in \mathcal{S}} \gamma_j \widehat{V}_{\cdot j}\right\|_2^2 (\Theta_{11} + \beta^2 \Theta_{22} - 2\beta \Theta_{12})$. Note that $T^{\beta} \mid W \sim N(0, \operatorname{Var}_1)$.*

Note that the difference between Lemma 10 from Lemma 9 is that $\mathcal{V}$ in Lemma 9 is replaced by $\mathcal{S}$ in Lemma 10. Using the same argument as Section 5.3, we establish Lemma 6.

## 5.6 Proof of Lemma 5

In the following, we present the proof of (68). The proof of (69) is similar and omitted here. By (65) and (66), on the event $\mathcal{B}$, we have

$$z_{\alpha/2}(1 - g(n)) \leq z_{\alpha/2} \sqrt{\frac{\widehat{\operatorname{Var}}(\widehat{\Sigma}_{12})}{\Theta_{22}^2 \operatorname{Var}_1 + \operatorname{Var}_2}} - \frac{R_1 + R_2 + R_3}{\sqrt{\operatorname{Var}_2 + \Theta_{22}^2 \operatorname{Var}_1}} \leq z_{\alpha/2}(1 + g(n)),$$

where $g(n) = C\left(\sqrt{s_{z1}} \frac{s \log p}{\sqrt{n}} + \frac{1}{\|\gamma\|_2} \frac{s_{z1} \log p}{\sqrt{n}}\right) + C \frac{1}{\sqrt{s_{z1} \log p}}$. Define the following events

$$\mathcal{F}_0 = \left\{\frac{M_1 + M_2}{\sqrt{\operatorname{Var}_2 + \Theta_{22}^2 \operatorname{Var}_1}} + B(\theta, W) \geq z_{\alpha/2} \sqrt{\frac{\widehat{\operatorname{Var}}(\widehat{\Sigma}_{12})}{\Theta_{22}^2 \operatorname{Var}_1 + \operatorname{Var}_2}} - \frac{R_1 + R_2 + R_3}{\sqrt{\operatorname{Var}_2 + \Theta_{22}^2 \operatorname{Var}_1}}\right\}$$

$$\mathcal{F}_1 = \left\{\frac{M_1 + M_2}{\sqrt{\operatorname{Var}_2 + \Theta_{22}^2 \operatorname{Var}_1}} + B(\theta, W) \geq z_{\alpha/2}(1 - g(n))\right\}$$

$$\mathcal{F}_2 = \left\{\frac{M_1 + M_2}{\sqrt{\operatorname{Var}_2 + \Theta_{22}^2 \operatorname{Var}_1}} + B(\theta, W) \geq z_{\alpha/2}(1 + g(n))\right\}$$

Note that $\mathcal{F}_2 \cap \mathcal{B} \subset \mathcal{F}_0 \cap \mathcal{B} \subset \mathcal{F}_1 \cap \mathcal{B}$. Hence, we have

$$\begin{aligned}
\left|\mathbf{P}\left(\mathcal{F}_{0}\right)-\left(1-\mathbf{E}_{W} \Phi\left(z_{\alpha / 2}-B(\theta, W)\right)\right)\right| & \leq \mathbf{P}\left(\mathcal{F}_{0} \cap \mathcal{B}^{c}\right) \\
& +\left|\mathbf{P}\left(\mathcal{F}_{0} \cap \mathcal{B}\right)-\left(1-\mathbf{E}_{W} \Phi\left(z_{\alpha / 2}-B(\theta, W)\right)\right)\right| \\
& \leq \mathbf{P}\left(\mathcal{F}_{0} \cap \mathcal{B}^{c}\right)+\max _{i=1,2}\left|\mathbf{P}\left(\mathcal{F}_{i} \cap \mathcal{B}\right)-\left(1-\mathbf{E}_{W} \Phi\left(z_{\alpha / 2}-B(\theta, W)\right)\right)\right|,
\end{aligned} \tag{102}$$

where the first inequality is an triangle inequality and the last inequality follows from the fact that $\mathcal{F}_2 \cap \mathcal{B} \subset \mathcal{F}_0 \cap \mathcal{B} \subset \mathcal{F}_1 \cap \mathcal{B}$. Note that

$$\begin{aligned}
& \left|\mathbf{P}\left(\mathcal{F}_{i} \cap \mathcal{B}\right)-\left(1-\mathbf{E}_{W} \Phi\left(z_{\alpha / 2}-B(\theta, W)\right)\right)\right| \\
& =\left|\mathbf{P}\left(\mathcal{F}_{i}\right)+\mathbf{P}(\mathcal{B})-\mathbf{P}\left(\mathcal{F}_{i} \cup \mathcal{B}\right)-\left(1-\Phi\left(z_{\alpha / 2}-B(\theta, W)\right)\right)\right| \\
& \leq\left|\mathbf{P}\left(\mathcal{F}_{i}\right)-\left(1-\mathbf{E}_{W} \Phi\left(z_{\alpha / 2}-B(\theta, W)\right)\right)\right|+\left|\mathbf{P}(\mathcal{B})-\mathbf{P}\left(\mathcal{F}_{i} \cup \mathcal{B}\right)\right|.
\end{aligned} \tag{103}$$

By (62), we have $\mathbf{P}\left(\mathcal{F}_i \cup \mathcal{B}\right) \to 1$ and $\mathbf{P}\left(\mathcal{F}_0 \cap \mathcal{B}^c\right) \to 0$ and hence it is sufficient to control the term $\max_{i=1,2}\left|\mathbf{P}\left(\mathcal{F}_i\right)-\left(1-\mathbf{E}_W \Phi\left(z_{\alpha/2}-B(\theta, W)\right)\right)\right|$. In the following, we will focus on the case $i=1$. The case for $i=2$ is similar. Note that

$$\begin{aligned}
\left|\mathbf{P}\left(\mathcal{F}_{1}\right)-\left(1-\mathbf{E}_{W} \Phi\left(z_{\alpha / 2}-B(\theta, W)\right)\right)\right| & \leq \mathbf{E}_{W}\left|\mathbf{P}\left(\mathcal{F}_{1} \mid W\right)-\left(1-\Phi\left(z_{\alpha / 2}-B(\theta, W)\right)\right)\right| \\
& =\mathbf{E}_{W}\left|\mathbf{P}\left(\mathcal{F}_{1} \mid W\right)-\left(1-\Phi\left(z_{\alpha / 2}-B(\theta, W)\right)\right)\right| \cdot \mathbf{1}_{\mathcal{B}_{8}} \\
& +\mathbf{E}_{W}\left|\mathbf{P}\left(\mathcal{F}_{1} \mid W\right)-\left(1-\Phi\left(z_{\alpha / 2}-B(\theta, W)\right)\right)\right| \cdot \mathbf{1}_{\mathcal{B}_{8}^{c}} \\
& \leq \mathbf{E}_{W}\left|\mathbf{P}\left(\mathcal{F}_{1} \mid W\right)-\left(1-\Phi\left(z_{\alpha / 2}-B(\theta, W)\right)\right)\right| \cdot \mathbf{1}_{\mathcal{B}_{8}}+2 \mathbf{P}\left(\mathcal{B}_{8}^{c}\right),
\end{aligned} \tag{104}$$

where the last inequality follows from $\left|\mathbf{P}\left(\mathcal{F}_1 \mid W\right)-\left(1-\Phi\left(z_{\alpha/2}-B(\theta, W)\right)\right)\right| \leq 2$. Starting form here, we will separate the proof into three cases,

Case a. $\|\gamma_{\mathcal{V}}\|_2 \ll \sqrt{n}$.

Case b. $\|\gamma_{\mathcal{V}}\|_2 \geq c\sqrt{n}$ and $\sqrt{n}\Sigma_{12} \to \Delta_1^*$;

Case c. $\|\gamma_{\mathcal{V}}\|_2 \geq c\sqrt{n}$ and $\sqrt{n}\Sigma_{12} \to \infty$;

<u>Case a.</u>

By (64), we have

$$\mathbf{P}(\mathcal{F}_1 \mid W) = \mathbf{E}_{M_2|W} \mathbf{P}(\mathcal{F}_1 \mid M_2, W)$$
$$= \mathbf{E}_{M_2|W} \left(1 - \Phi\left(\frac{z_{\alpha/2}(1-g(n)) + B(\theta, W)\sqrt{\text{Var}_2 + \Theta_{22}^2 \text{Var}_1} - M_2}{\sqrt{\Theta_{22}^2 \text{Var}_1}}\right)\right).$$

By the proposition [KMT] in Mason et al. (2012), conditioning on $W$, there exist $\bar{M}_2$ having the same distribution as $M_2$ and $\widetilde{M}_2 \sim N(0, \text{Var}_2)$ such that on the event $\mathcal{F}_4 = \left\{\left|\frac{\widetilde{M}_2}{\sqrt{\text{Var}_2}}\right| \leq C_1 \sqrt{n}\right\}$, we have $\left|\bar{M}_2 - \widetilde{M}_2\right| \leq C_2 \left(\frac{\bar{M}_2^2}{\sqrt{n \text{Var}_2}} + \sqrt{\frac{\text{Var}_2}{n}}\right)$. Hence we have

$$\mathbf{P}(\mathcal{F}_1 \mid W) = \mathbf{E}_{M_2|W} \left(1 - \Phi\left(\frac{z_{\alpha/2}(1-g(n)) + B(\theta, W)\sqrt{\text{Var}_2 + \Theta_{22}^2 \text{Var}_1} - M_2}{\sqrt{\Theta_{22}^2 \text{Var}_1}}\right)\right)$$
$$= \mathbf{E}_{\bar{M}_2|W} \left(1 - \Phi\left(\frac{z_{\alpha/2}(1-g(n)) + B(\theta, W)\sqrt{\text{Var}_2 + \Theta_{22}^2 \text{Var}_1} - \bar{M}_2}{\sqrt{\Theta_{22}^2 \text{Var}_1}}\right)\right)$$
$$= \mathbf{E}_{\widetilde{M}_2|W} \left(1 - \Phi\left(\frac{z_{\alpha/2}(1-g(n)) + B(\theta, W)\sqrt{\text{Var}_2 + \Theta_{22}^2 \text{Var}_1} - \widetilde{M}_2}{\sqrt{\Theta_{22}^2 \text{Var}_1}}\right)\right) + \mathbf{E}_{\bar{M}_2, \widetilde{M}_2|W} g_1(n),$$
(105)

where

$$g_1(n) = \Phi\left(\frac{z_{\alpha/2}(1-g(n)) + B(\theta, W)\sqrt{\text{Var}_2 + \Theta_{22}^2 \text{Var}_1} - \widetilde{M}_2}{\sqrt{\Theta_{22}^2 \text{Var}_1}}\right) - \Phi\left(\frac{z_{\alpha/2}(1-g(n)) + B(\theta, W)\sqrt{\text{Var}_2 + \Theta_{22}^2 \text{Var}_1} - \bar{M}_2}{\sqrt{\Theta_{22}^2 \text{Var}_1}}\right)$$

Since

$$1 - \Phi\left(z_{\alpha/2} - B(\theta, W)\right) = \mathbf{E}_{\widetilde{M}_2|W}\left(1 - \Phi\left(\frac{z_{\alpha/2}(1-g(n)) + B(\theta, W)\sqrt{\text{Var}_2 + \Theta_{22}^2 \text{Var}_1} - \widetilde{M}_2}{\sqrt{\Theta_{22}^2 \text{Var}_1}}\right)\right),$$

it follows from (105) that

$$\left|\mathbf{P}(\mathcal{F}_1 \mid W) - \left(1 - \Phi\left(z_{\alpha/2} - B(\theta, W)\right)\right)\right| = \left|\mathbf{E}_{\bar{M}_2, \widetilde{M}_2|W} g_1(n)\right| \qquad (106)$$

Note that

$$\left|\mathbf{E}_{\bar{M}_2, \widetilde{M}_2|W} g_1(n)\right| \leq \mathbf{E}_{\bar{M}_2, \widetilde{M}_2|W} |g_1(n)| \cdot \mathbf{1}_{\mathcal{F}_4} + \mathbf{E}_{\bar{M}_2, \widetilde{M}_2|W} |g_1(n)| \cdot \mathbf{1}_{\mathcal{F}_4^c} \leq \mathbf{E}_{\bar{M}_2, \widetilde{M}_2|W} |g_1(n)| \cdot \mathbf{1}_{\mathcal{F}_4} + 2\mathbf{P}(\mathcal{F}_4^c).$$

and

$$\mathbf{E}_{\bar{M}_2,\widetilde{M}_2|W}|g_1(n)|\cdot \mathbf{1}_{\mathcal{F}_4} \le 2\mathbf{E}_{\bar{M}_2,\widetilde{M}_2|W}\frac{\left|\bar{M}_2-\widetilde{M}_2\right|}{\sqrt{\Theta_{22}^2\mathrm{Var}_1}} \le 2C_2\frac{1}{\sqrt{\Theta_{22}^2\mathrm{Var}_1}}\mathbf{E}_{\bar{M}_2,\widetilde{M}_2|W}\left(\frac{\bar{M}_2^2}{\sqrt{n\mathrm{Var}_2}}+\sqrt{\frac{\mathrm{Var}_2}{n}}\right).$$

Combined with (106), we have

$$\left|\mathbf{P}\left(\mathcal{F}_1 \mid W\right) - \left(1-\Phi\left(z_{\alpha/2}-B\left(\theta,W\right)\right)\right)\right| \le 2C_2\frac{1}{\sqrt{\Theta_{22}^2\mathrm{Var}_1}}\mathbf{E}_{\bar{M}_2,\widetilde{M}_2|W}\left(\frac{\bar{M}_2^2}{\sqrt{n\mathrm{Var}_2}}+\sqrt{\frac{\mathrm{Var}_2}{n}}\right)+2\mathbf{P}\left(\mathcal{F}_4^c\right)$$

$$= 4C_2\frac{1}{\sqrt{\Theta_{22}^2\mathrm{Var}_1}}\sqrt{\frac{\mathrm{Var}_2}{n}}+2\mathbf{P}\left(\mathcal{F}_4^c\right),$$

where the last equality follows from the fact that $\mathbf{E}\bar{M}_2^2 = \mathrm{Var}_2$. Combined with (104), we have

$$\left|\mathbf{P}\left(\mathcal{F}_1\right)-\left(1-\mathbf{E}_W\Phi\left(z_{\alpha/2}-B\left(\theta,W\right)\right)\right)\right| \le \mathbf{E}_W\left(4C_2\frac{1}{\sqrt{\Theta_{22}^2\mathrm{Var}_1}}\sqrt{\frac{\mathrm{Var}_2}{n}}\cdot\mathbf{1}_{\mathcal{B}_8}\right)$$
$$+2\mathbf{P}\left(\mathcal{F}_4^c\right)+2\mathbf{P}\left(\mathcal{B}_8^c\right) \le C\frac{\|\gamma_\mathcal{V}\|_2}{\sqrt{\Theta_{22}^2}}\sqrt{\frac{\mathrm{Var}_2}{n}}+2\mathbf{P}\left(\mathcal{F}_4^c\right)+2\mathbf{P}\left(\mathcal{B}_8^c\right),$$

where the last inequality follows from the definition of $\mathcal{B}_8$. Under the assumption $\|\gamma_\mathcal{V}\|_2 \ll \sqrt{n}$, we show that $\left|\mathbf{P}\left(\mathcal{F}_1\right)-\left(1-\mathbf{E}_W\Phi\left(z_{\alpha/2}-B\left(\theta,W\right)\right)\right)\right| \to 0$. Combined with (102) and (103), we establish (68).

Case b.

Under the assumption $\|\gamma_\mathcal{V}\|_2 \ge c\sqrt{n}$, on the event $\mathcal{B}_8$, we have $\mathrm{Var}_1 \to 0$ and

$$M_1 \overset{p}{\to} 0, \quad M_2 \overset{d}{\to} N(0,\mathrm{Var}_2), \quad B(\theta,W) \to \frac{\Delta_1^*}{\sqrt{\mathrm{Var}_2}}.$$

By the bounded convergence theorem, we establish

$$\mathbf{E}_W\left|\mathbf{P}\left(\mathcal{F}_1 \mid W\right)-\left(1-\Phi\left(z_{\alpha/2}-\frac{\Delta_1^*}{\sqrt{\mathrm{Var}_2}}\right)\right)\right|\cdot\mathbf{1}_{\mathcal{B}_8} \to 0$$
$$\mathbf{E}_W\left|\Phi\left(z_{\alpha/2}-\frac{\Delta_1^*}{\sqrt{\mathrm{Var}_2}}\right)-\Phi\left(z_{\alpha/2}-\frac{\Delta_1^*}{\sqrt{\mathrm{Var}_2}}\right)\right|\cdot\mathbf{1}_{\mathcal{B}_8} \to 0.$$

By applying triangle inequality, we establish (68).

Case c.

This case is similar to Case b. The only difference is $B(\theta, W) \to \infty$ and on the event $\mathcal{B}_8$, we have

$$\mathbf{P}\left(\mathcal{F}_1 \mid W\right) \cdot \mathbf{1}_{\mathcal{B}_8} \to 1, \quad \Phi\left(z_{\alpha/2} - \frac{\Delta_1^*}{\sqrt{\mathrm{Var}_2}}\right) \cdot \mathbf{1}_{\mathcal{B}_8} \to 0.$$